\documentclass[10pt]{article}
\usepackage{xcolor}
\usepackage{amssymb}
\usepackage{cite}

\newlength{\minitwocolumn}
\font\teneufm=eufm10
\font\seveneufm=eufm7
\font\fiveeufm=eufm5
\newfam\eufmfam
\textfont\eufmfam=\teneufm
\scriptfont\eufmfam=\seveneufm
\scriptscriptfont\eufmfam=\fiveeufm

\renewcommand{\theequation}{\arabic{equation}}
\newtheorem{thm}{Theorem}[section]

\newtheorem{dfn}[thm]{Definition}
\newtheorem{prop}[thm]{Proposition}
\newtheorem{lem}[thm]{Lemma}

\usepackage{amsmath,amssymb}
\title{
\large{\bf
Quadratic relations of the deformed $W$-superalgebra 
${\cal W}_{q, t}\bigl(A(M,N)\bigr)$}}
\author{Takeo KOJIMA}
\begin{document}

\maketitle

\begin{center}
{\it Department of Mathematics and Physics, Faculty of Engineering, Yamagata University,\\
Jonan 4-chome 3-16, Yonezawa 992-8510, JAPAN}
\end{center}

\begin{abstract}
We find the free field construction of the basic $W$-current
and screening currents for the deformed $W$-superalgebra
${\cal W}_{q,t}\bigl(A(M,N)\bigr)$ associated with Lie superalgebra of type $A(M,N)$.
Using this free field construction,
we introduce the higher $W$-currents and obtain a closed set
of quadratic relations among them.
These relations are independent of the choice of Dynkin-diagrams
for the Lie superalgebra $A(M,N)$, though the screening currents 
are not. This allows us to define ${\cal W}_{q,t}\bigl(A(M,N)\bigr)$ by generators and relations.
\end{abstract}

\section{Introduction}
\label{Section1}

The deformed $W$-algebra ${\cal W}_{q,t}\bigl(\mathfrak{g}\bigr)$ is a two parameter deformation
of the classical $W$-algebra ${\cal W}(\mathfrak{g})$.
Shiraishi et al. \cite{Shiraishi-Kubo-Awata-Odake} obtained a free field
construction of the deformed Virasoro algebra
${\cal W}_{q,t}\bigl(\mathfrak{sl}(2)\bigr)$, which is a one-parameter deformation of the Virasoro algebra,
to construct a deformation of the correspondence between conformal field theory and the Calogero-Sutherland model.
The theory of the deformed $W$-algebras 
${\cal W}_{q,t}(\mathfrak{g})$
has been developed in papers
\cite{
Awata-Kubo-Odake-Shiraishi, Feigin-Frenkel, Brazhnikov-Lukyanov, Hara-Jimbo-Konno-Odake-Shiraishi, Frenkel-Reshetikhin, Sevostyanov, Odake, Feigin-Jimbo-Mukhin-Vilkoviskiy, Ding-Feigin, Kojima}.
However, in comparison with the conformal case,
the theory of the deformed $W$-algebra 
is still not fully developed and understood.
For that matter it is worthwhile to concretely 
construct ${\cal W}_{q,t}({\mathfrak g})$ in each case.
This paper is a continuation
of the paper \cite{Kojima} for ${\cal W}_{q,t}\bigl(A(1,0)\bigr)$.
The purpose of this paper is to generalize the result of case 
$A(1,0)$ to $A(M,N)$.

We follow the method of 
\cite{Ding-Feigin}, where a free field construction is found for 
the deformed $\mathcal{W}_{q,t}\bigl(\mathfrak{sl}(3)\bigr)$ and  $\mathcal{W}_{q,t}\bigl(A(1,0)\bigr)$. 
Starting from a $W$ current given as a sum of three vertex operators 
\begin{eqnarray}
T_1(z)=\Lambda_1(z)+\Lambda_2(z)+\Lambda_3(z)\,,
\nonumber
\end{eqnarray}
and two screening currents $S_j(z)$ given by a vertex operator, 
the authors of \cite{Ding-Feigin} determined them simultaneously by demanding that $T_1(z)$ and $S_j(w)$ commute up to a total difference.
Higher currents $T_i(z)$ are defined inductively by the fusion relation 
\begin{eqnarray}
\mathop{\mathrm{Res}}_{w=x^i z}T_1(w)T_{i-1}(z)=c_i T_i(x^{i-1}z)
\nonumber
\end{eqnarray}
with appropriate constants $x$ and $c_i$. In the case of $\mathcal{W}_{q,t}\bigl(\mathfrak{sl}(3)\bigr)$,  
it is known that they truncate, {\it i.e.} $T_3(z)=1$ and $T_i(z)=0$ ($i\ge4$),
and that $T_1(z)$ and $T_2(z)$ satisfy the quadratic relations 
\cite{Awata-Kubo-Odake-Shiraishi, Feigin-Frenkel}
\begin{eqnarray}
&&
f_{1,1}\left(\frac{z_2}{z_1}\right)T_1(z_1)T_1(z_2)-
f_{1,1}\left(\frac{z_1}{z_2}\right)T_1(z_2)T_1(z_1)
=c\left(
\delta\left(\frac{x^{-2}z_2}{z_1}\right)T_2(x^{-1}z_2)
-\delta\left(\frac{x^2z_2}{z_1}\right)T_2(x z_2)
\right),\nonumber\\
&&
f_{1,2}\left(\frac{z_2}{z_1}\right)T_1(z_1)T_2(z_2)-
f_{2,1}\left(\frac{z_1}{z_2}\right)T_2(z_2)T_1(z_1)
=c\left(
\delta\left(\frac{x^{-3}z_2}{z_1}\right)
-\delta\left(\frac{x^3z_2}{z_1}\right)
\right),
\nonumber
\\
&&
f_{2,2}\left(\frac{z_2}{z_1}\right)T_2(z_1)T_2(z_2)-
f_{2,2}\left(\frac{z_1}{z_2}\right)T_2(z_2)T_2(z_1)
=c\left(
\delta\left(\frac{x^{-2}z_2}{z_1}\right)T_1(x^{-1}z_2)
-\delta\left(\frac{x^2z_2}{z_1}\right)T_1(x z_2)
\right)\nonumber
\end{eqnarray}
with appropriate constants $x, c$, and functions $f_{i,j}(z)$.
In the case of $\mathcal{W}_{q,t}\bigl(A(1,0)\bigr)$, 
it was shown in \cite{Kojima} that
such truncation for $T_i(z)$ does not
take place and that an infinite number of quadratic relations are satisfied by infinite number of $T_i(z)$'s.
In the present paper, we extend this result to general $A(M,N)$.

Following the method of \cite{Ding-Feigin},
we construct the basic $W$ current $T_1(z)$ together with the screening currents $S_j(w)$
for  $\mathcal{W}_{q,t}\bigl(A(M,N)\bigr)$
(See (\ref{def:T1(z)}) and (\ref{def:screening})).
We introduce the higher $W$-currents $T_i(z)$
(See (\ref{def:Ti(z)})) and obtain a closed set
of quadratic relations among them (See (\ref{thm:quadratic})).
We show further that these relations are independent of the choice 
of Dynkin-diagrams
for the superalgebra $A(M,N)$, though the screening currents 
are not. This allows us to define ${\cal W}_{q,t}\bigl(A(M,N)\bigr)$ by generators and relations.

The text is organized as follows.
In Section \ref{Section2},  
we prepare the notation
and formulate the problem.
In Section \ref{Section3},
we give a free field construction of
the basic $W$-current $T_1(z)$ and the screening currents 
$S_j(w)$ for the deformed $W$-algebra 
${\cal W}_{q,t}\bigl(A(M,N)\bigr)$. 
In Section \ref{Section4}, we introduce higher $W$-currents 
$T_i(z)$ 
and present a closed set of quadratic relations among them. 
We show
that these quadratic relations are independ of the choice of the Dynkin-diagram for the superalgebra $A(M,N)$.
We also obtain the $q$-Poisson algebra in the classical limit. 
Section \ref{Section5} is devoted to conclusion and discussion.

\section{Preliminaries}
\label{Section2}

In this section we prepare the notation and formulate the problem.
Throughout this paper, we fix a real number $r>1$ 
and a complex number $x$ with $0<|x|<1$.

\subsection{Notation}

In this section we use complex numbers 
$a$, $w$ $(w\neq 0)$, 
$q$ ($q \neq 0,\pm1$), 
and $p$ with $|p|<1$.
For any integer $n$, define $q$-integer
\begin{eqnarray}
[n]_q=\frac{q^n-q^{-n}}{q-q^{-1}}.
\nonumber
\end{eqnarray}
We use symbols for infinite products
\begin{eqnarray}
(a;p)_\infty=\prod_{k=0}^\infty (1-a p^k),~~~
(a_1,a_2, \ldots, a_N; p)_\infty=
\prod_{i=1}^N (a_i; p)_\infty
\nonumber
\end{eqnarray}
for complex numbers $a_1, a_2, \ldots, a_N$.
The following standard formulae are useful.
\begin{eqnarray}
\exp\left(-\sum_{m=1}^\infty \frac{1}{m}a^m \right)=1-a,~~~\exp\left(-\sum_{m=1}^\infty \frac{1}{m}\frac{a^m}{1-p^m}\right)=(a;p)_\infty.\nonumber
\end{eqnarray}
We use the elliptic theta function $\Theta_p(w)$ and the compact notation $\Theta_p(w_1,w_2, \ldots, w_N)$ as
\begin{eqnarray}
\Theta_p(w)=(p, w, p w^{-1};p)_\infty,~~~
\Theta_p (w_1, w_2, \ldots, w_N)=\prod_{i=1}^N \Theta_p(w_i)
\nonumber
\end{eqnarray}
for complex numbers $w_1, w_2, \ldots, w_N \neq 0$.
Define 
$\delta(z)$ by the formal series
\begin{eqnarray}
\delta(z)=\sum_{m \in {\mathbf Z}}z^m.
\nonumber
\end{eqnarray}

\subsection{Dynkin-diagram of $A(M,N)$}

In this section we introduce Dynkin-diagrams of 
the Lie superalgebra
$A(M,N)$.
We fix integers $M, N$ $(M+N \geq 1, M, N=0,1,2,\ldots)$.
We set $L=M+N+1$.
Let $\varepsilon_1, \varepsilon_2, \ldots, \varepsilon_{M+1}$
and $\delta_1, \delta_2, \ldots, \delta_{N+1}$
be a basis of $\mathbf{R}^{L+1}$ with an inner product
$(~,~)$ such that
\begin{eqnarray}
&&
(\varepsilon_i, \varepsilon_j)=\delta_{i,j}~~~(1\leq i, j \leq M+1),~~~
(\delta_i, \delta_j)=-\delta_{i,j}~~~(1\leq i, j \leq N+1),
\nonumber\\
&&
(\varepsilon_i, \delta_j)=(\delta_j, \varepsilon_i)=0~~~
(1\leq i \leq M+1, 1\leq j \leq N+1).
\nonumber
\end{eqnarray}
The standard fundamental system $\Pi^{st}$
for the Lie superalgebra $A(M,N)$ is given as
\begin{eqnarray}
&&\Pi^{st}=
\{\alpha_i=\varepsilon_i-\varepsilon_{i+1},
\alpha_{M+1}=\varepsilon_{M+1}-\delta_1,
\alpha_{M+1+j}=\delta_j-\delta_{j+1}
|1\leq i \leq M, 1\leq j \leq N\}.\nonumber
\end{eqnarray}
The standard Dynkin-diagram $\Phi^{st}$ for 
the Lie superalgebra 
$A(M, N)$ is given as\\
\begin{center}
\input{Dynkin1.tex}
\end{center}
~\\
Here a circle represents an even simple root and
a crossed circle represents an odd isotropic simple root.

There is an indeterminacy in how to choose 
Dynkin-diagram for the Lie superalgebra
$A(M,N)$, which is brought by 
fundamental reflections $r_{\alpha_i}$.
For the fundamental system $\Pi$,
the fundamental reflection $r_{\alpha_i}$ $(\alpha_i \in \Pi)$
satisfies
\begin{eqnarray}
r_{\alpha_i}(\alpha_j)=
\left\{\begin{array}{cc}
-\alpha_i&{\rm if}~~~j=i,\\
\alpha_i+\alpha_j&{\rm if}~~~j \neq i,~ (\alpha_i, \alpha_j)\neq 0,
\\
\alpha_j&{\rm if}~~~j \neq i,~ (\alpha_i, \alpha_j)=0.
\end{array}
\right.
\nonumber
\end{eqnarray}
For an odd isotropic root $\alpha_i$,
we call the fundamental reflection $r_{\alpha_i}$
odd reflection.
For an even root $\alpha_i$,
we call the fundamental reflection $r_{\alpha_i}$
real reflection.
The Dynkin-diagram transformed by $r_{\alpha_i}$ is
represented as $r_{\alpha_i}(\Phi)$.
Real reflections don't change Dynkin-diagram.
We illustrate the notion of odd reflections as follows.\\

\begin{center}
\input{Dynkin2.tex}~~~~~
\input{Arrow1.tex}~~~
\input{Dynkin3.tex}
\\~\\
\input{Dynkin4.tex}~~~~
\input{Arrow1.tex}~~~
\input{Dynkin5.tex}
\\~
\\
\input{Dynkin6.tex}~~~~~~
\input{Arrow2.tex}~~~
\input{Dynkin7.tex}
\\~
\\
\input{Dynkin8.tex}~~~~~~
\input{Arrow3.tex}~~~
\input{Dynkin9.tex}
\\~
\\
\end{center}
{\bf Example}~~$A(1,0)$ and $A(0,1)$\\
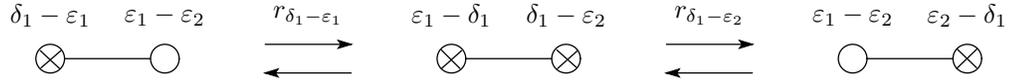
\begin{center}
\input{Dynkin10.tex}
\end{center}~\\
Here $\Pi_1=\{\delta_1-\varepsilon_1, \varepsilon_1-\varepsilon_2\}$ and $\Pi_2=\{\varepsilon_1-\delta_1, \delta_1-\varepsilon_2\}$ are the other fundamental systems.
\\~\\
{\bf Example}~~$A(1,1)$\\
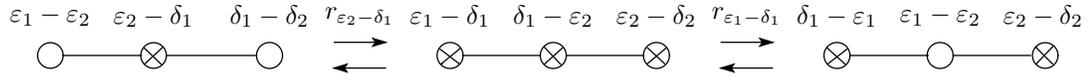
\begin{center}
\input{Dynkin11.tex}
\end{center}
~\\
Here $\Pi_1=\{\varepsilon_1-\delta_1, 
\delta_1-\varepsilon_2, \varepsilon_2-\delta_2\}$ 
and 
$\Pi_2=\{\delta_1-\varepsilon_1, \varepsilon_1-\varepsilon_2,
\varepsilon_2-\delta_2\}$ are the other fundamental systems.
~\\
\\
{\bf Example}~~$A(2,0)$ and $A(0,2)$\\
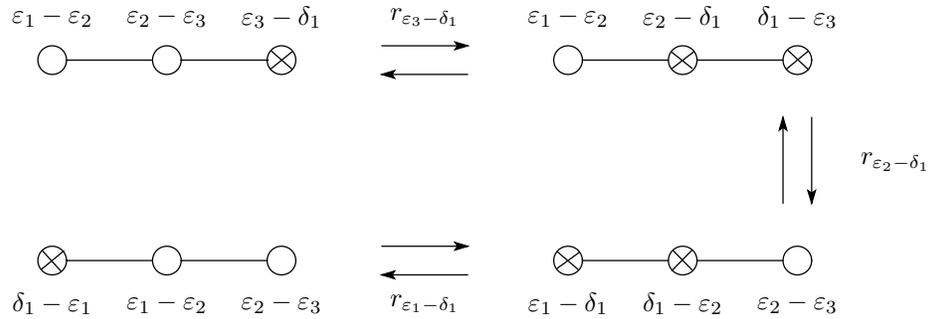
\begin{center}
\input{Dynkin12.tex}
\end{center}
~\\
Here $\Pi_1=\{\varepsilon_1-\varepsilon_2, 
\varepsilon_2-\delta_1, \delta_1-\varepsilon_3\}$,
$\Pi_2=\{\varepsilon_1-\delta_1, \delta_1-\varepsilon_2,
\varepsilon_2-\varepsilon_3\}$,
and
$\Pi_3=\{\delta_1-\varepsilon_1, \varepsilon_1-\varepsilon_2,
\varepsilon_2-\varepsilon_3\}$
are the other fundamental systems.

\subsection{Ding-Feigin's construction}

We introduce the Heisenberg algebra with generators
$a_i(m)$, $Q_i$ $(m \in {\mathbf Z}, 1\leq i \leq L)$ satisfying
\begin{eqnarray}
&&
[a_i(m), a_j(n)]=\frac{1}{m}A_{i,j}(m)\delta_{m+n,0}~~
(m, n \neq 0, 1\leq i,j \leq L),
\nonumber\\
&&
[a_i(0),Q_j]=A_{i,j}(0)~~(1\leq i,j \leq L).
\nonumber
\end{eqnarray}
The remaining commutators vanish.
We impose the following conditions on the parameters $A_{i,j}(m)\in\mathbf{C}$:
\begin{eqnarray}
&&
A_{i,i}(m)=1~(m \neq 0, 1\leq i \leq L),~~~A_{i,j}(m)=A_{j,i}(-m)~(m \in \mathbf{Z}, 1\leq i \neq  j \leq L),
\nonumber\\
&&
\det\left(\left(A_{i,j}(m)\right)_{i,j=1}^L\right)\neq 0~~(m \in \mathbf{Z}).
\nonumber
\end{eqnarray}
We use the normal ordering symbol $:~:$ that satisfies
\begin{eqnarray}
&&
:a_i(m)a_j(n):=\left\{\begin{array}{cc}
a_i(m)a_j(n)&(m<0),\\
a_j(n)a_i(m)&(m \geq 0)
\end{array}
\right.~~~(m,n \in {\mathbf Z}, 1\leq i,j \leq L),
\nonumber
\\
&&
:a_i(0)Q_j:=:Q_j a_i(0):=Q_j a_i(0)~~~(1 \leq i,j \leq L).
\nonumber
\end{eqnarray}

Next, we work on Fock space of the free field.
Let $T_1(z)$ be a sum of vertex operators
\begin{align}
&T_1(z)=g_1 \Lambda_1(z)+g_2 \Lambda_2(z)+\cdots+g_{L+1} \Lambda_{L+1}(z),
\label{def:T1(z)}\\
&{\Lambda}_i(z)=e^{\sum_{j=1}^L \lambda_{i,j}(0)a_j(0)}
:\exp\left(\sum_{j=1}^L \sum_{m \neq 0}\lambda_{i,j}(m)a_j(m)z^{-m}\right):~~(1\leq i \leq L+1).
\nonumber
\end{align}
We call  $T_1(z)$ the basic $W$-current.
We introduce the screening currents 
$S_j(w)$ $(1 \leq j \leq L)$ as
\begin{eqnarray}
S_j(w)=w^{\frac{1}{2}A_{j,j}(0)}e^{Q_j}w^{a_j(0)}
:\exp\left(\sum_{m \neq 0}s_j(m)a_j(m)w^{-m}\right):
~~(1\leq j \leq L).
\label{def:screening}
\end{eqnarray}
The parameters
$A_{i, j}(m)$, $\lambda_{i, j}(m), s_j(m)$ and $g_i$ are to be determined through the construction given below.

Quite generally, given two vertex  operators $V(z)$, $W(w)$, 
their product has the form 
\begin{eqnarray}
V(z)W(w)=\varphi_{V,W}\left(z,w\right):V(z)W(w):~~~(|z|\gg |w|)
\nonumber
\end{eqnarray}
with some formal power series $\varphi_{V,W}(z,w)\in\mathbf{C}[[w/z]]$. 
The vertex operators $V(z)$ and $W(w)$ are said to be mutually local if the following two conditions hold.
\begin{align}
&\mathrm{(i)}~~\varphi_{V,W}(z,w)~{\rm and}~\varphi_{W,V}(w,z)~{\rm converge~ to~rational~functions},~~~\nonumber\\
&\mathrm{(ii)}~~\varphi_{V,W}(z,w)=\varphi_{W,V}(w,z).\nonumber
\end{align}

Under this setting, we are going to determine
the $W$-current $T_1(z)$ 
and the screening currents $S_j(w)$ that
satisfy the following
mutual locality (\ref{assume:mutual-locality}), commutativity (\ref{assume:q-difference}),
and symmetry (\ref{assume:varphiSS}).
\\
\\
{\bf Mutual Locality}~~~$\Lambda_i(z)$ $(1\leq i \leq L+1)$ and 
$S_j(w)$ $(1\leq j \leq L)$ are mutually local and the operator product expansion
have at most one pole and one zero.
\begin{eqnarray}
\varphi_{\Lambda_i, S_j}(z,w)=\varphi_{S_j, \Lambda_i}(w,z)
={\displaystyle \frac{w-\frac{z}{p_{i,j}}}{w-\frac{z}{q_{i,j}}}}~~~(1\leq i \leq L+1, 1\leq j \leq L).
\label{assume:mutual-locality}
\end{eqnarray}
We allow the possibility $p_{i,j}=q_{i,j}$,  in which case 
$\Lambda_i(z)S_j(w)=S_j(w)\Lambda_i(z)=:\Lambda_i(z)S_j(w):$.

~\\
{\bf Commutativity}~~~$T_1(z)$ commutes with $S_j(w)$ $(1\leq j \leq L)$
up to a total difference
\begin{eqnarray}
[T_1(z),S_j(w)]=B_j(z)
\left(\delta\left(\frac{q_{j,j}w}{z}\right)-\delta\left(\frac{q_{j+1,j}w}{z}\right)\right)~~~(1\leq j \leq L)\,,
\label{assume:q-difference}
\end{eqnarray}
with some currents $B_j(z)$ $(1\leq j \leq L)$.

~\\
{\bf Symmetry}~~~For $\widetilde{S}_j(w)=e^{-Q_j}S_j(w)$ $(1\leq j \leq L)$, we impose
\begin{eqnarray}
&&
\varphi_{\widetilde{S}_k, \widetilde{S}_{l}}(w,z)=
\varphi_{\widetilde{S}_{l}, \widetilde{S}_k}(w,z)~~~
(1\leq k, l \leq L),
\nonumber\\
&&
\varphi_{\widetilde{S}_k, \widetilde{S}_l}(w,z)=1~~~
(|k-l|\geq 2, 1\leq k, l \leq L).
\label{assume:varphiSS}
\end{eqnarray}
For simplicity, we impose further the following conditions.
\begin{align}
&q_{i,j}~~(1\leq i \leq L+1, 1\leq j \leq L)~~{\rm are~distinct},
\label{assume:p,q0}\\
&\left|\frac{q_{j+1,j}}{q_{j,j}}\right| \neq 1~~~(1\leq j \leq L),
~~~-1<A_{k,k+1}(0)<0~~~(1\leq k \leq L-1).
\label{assume:p,q1}
\end{align}

It can be seen from elementary consideration that
there are three kinds of freedom in choosing parameters.
\\
{(i)}~Rearranging indices
\begin{eqnarray}
\Lambda_{i}(z) \mapsto \Lambda_{i'}(z),~~~S_j(w) \mapsto S_{j'}(w).
\label{def:arrange-suffix}
\end{eqnarray}
{(ii)}~Scaling variables: $\Lambda_i(z) \mapsto \Lambda_i(sz)$, i.e.
\begin{eqnarray}
\lambda_{i,j}(m) \mapsto s^m \lambda_{i,j}(m)~~~
q_{i,j} \mapsto s q_{i,j},~~~p_{i,j} \mapsto s p_{i,j}~~~(m\neq 0, 1\leq i \leq L+1, 1\leq j \leq L).
\label{def:shift}
\end{eqnarray}
{(iii)}~Scaling free fields: The free field can be rescaled as 
$a_j(m) \mapsto \alpha_j(m)^{-1} a_j(m)$ $(m>0, 1\leq j \leq L)$, i.e.
\begin{eqnarray}
\begin{array}{cc}
\begin{array}{c}
a_j(m) \mapsto \alpha_j(m)^{-1} a_j(m),~~s_j(m) \mapsto \alpha_j(m) s_j(m),\\
\lambda_{i,j}(m) \mapsto \lambda_{i,j}(m)\alpha_j(m)
\end{array}&(m \neq 0, 1\leq i \leq L+1, 1 \leq j \leq L),
\\
A_{i,j}(m) \mapsto \alpha_i(m)^{-1} A_{i,j}(m) \alpha_j(m)
&(m \neq 0, 1\leq i, j \leq L),
\end{array}
\label{def:scaling-boson}
\end{eqnarray}
where we set 
$\alpha_j(m)\neq 0$ $(m>0, 1\leq j \leq L)$ and $\alpha_j(-m)=\alpha_j(m)^{-1}$ $(m>0, 1\leq j \leq L)$.

\section{Free field construction}
\label{Section3}

In this section we give a free field construction of
the basic $W$-current and the screening currents for
${\cal W}_{q,t}\bigl(A(M,N)\bigr)$.

\subsection{Free field construction}

In Ding-Feigin's construction\cite{Ding-Feigin}, 
there are $2^L$ cases to be considered separately
according to values of $A_{j,j}(0)$ $(1\leq j \leq L)$.
We fix a pair of integers $j_1, j_2, \ldots , j_K$ $(1 \leq K \leq L)$
satisfying $1\leq j_1 <j_2<\cdots <j_K \leq L$.
Hereafter, 
we study the case the following conditions for
$A_{j,j}(0)$ $(1\leq j \leq L)$
are satisfied.
\begin{eqnarray}
~A_{j,j}(0)=1~~~{\rm if}~~~j=j_1, j_2, \ldots, j_K,~~
A_{j,j}(0)\neq 1~~~{\rm if}~~~j\neq j_1, j_2,\ldots, j_K.
\nonumber
\end{eqnarray}

First, we prepare the parameters $A_{i,j}(0)$ to give the free field construction.
We have already introduced $L\times L$ symmetric matrix $(A_{i,j}(0))_{i,j=1}^L$ as parameters of
the Heisenberg algebra.
To write $p_{i,j}$, $q_{i,j}$, $A_{i,j}(m)$, $s_j(m)$, 
and $\lambda_{i,j}(m)$ explicitly,
it is convenient to introduce $(L+1)\times (L+1)$ symmetric matrix $(A_{i,j}(0))_{i,j=0}^L$ 
uniquely extended from $(A_{i,j}(0))_{i,j=1}^{L}$ as follows.
\begin{eqnarray}
&&
A_{0,1}(0)=\left\{
\begin{array}{cc}
A_{1,2}(0)&{\rm if}~~j_1 \neq 1,\\
-1-A_{1,2}(0)&{\rm if}~~j_1=1,
\end{array}
\right.
~A_{0,L}(0)=\left\{
\begin{array}{cc}
A_{L, L-1}(0)&{\rm if}~~j_K \neq L,\\
-1-A_{L, L-1}(0)&{\rm if}~~j_K=L,
\end{array}
\right.
\nonumber
\\
&&
A_{0,0}(0)=\left\{\begin{array}{cc}
-2A_{0,L}(0)&{\rm if}~~K={\rm even},
\\
1&{\rm if}~~K={\rm odd},
\end{array}\right.
~~
A_{0,i}(0)=0~~(i \neq 0, 1, L).
\label{eqn:A(0)-0}
\end{eqnarray}
The extended matrix $\left(A_{i,j}(0)\right)_{i,j=0}^L$ are
explicitly written by $\beta=A_{1,2}(0)$
as follows (See Lemma \ref{lem:3-10}).
\begin{eqnarray}
&&
A_{i,i}(0)=\left\{
\begin{array}{cc}
1&{\rm if}~~i \in \widehat{J},\\
-2\beta&{\rm if}~~
i \notin \widehat{J},~i \in \widehat{I}(\beta),\\
2(1+\beta)&{\rm if}~~
i \notin \widehat{J},~i \in \widehat{I}(-1-\beta)
\end{array}
\right.~~(1 \leq i \leq L+1),
\nonumber
\\
&&
A_{j-1,j}(0)=A_{j,j-1}(0)=\left\{\begin{array}{cc}
\beta&{\rm if}~~j \in \widehat{I}(\beta),
\\
-1-\beta&{\rm if}~~j \in \widehat{I}(-1-\beta)
\end{array}\right.~~(1\leq j \leq L+1),
\nonumber
\\
&&
A_{k,l}(0)=A_{l,k}(0)=0~~\left(
|k-l|\geq 2, 1\leq k, l \leq L~~
{\rm or}~~k=0,~l \neq 0, 1, L\right).
\label{eqn:A(0)-1}
\end{eqnarray}
Here we set
\begin{eqnarray}
&&
\widehat{J}=\left\{\begin{array}{cc}
\{j_1, j_2, \ldots, j_K\}&{\rm if}~~K={\rm even},\\
\{j_1, j_2, \ldots, j_K, L+1\}&{\rm if}~~K={\rm odd},
\end{array}
\right.~~~
\widehat{I}(\delta)=
\{1\leq j \leq L+1|A_{j-1,j}(0)=\delta\}.\nonumber
\end{eqnarray}
We understand subscripts of $A_{i,j}(0)$ with mod.$L+1$, i.e.
$A_{0, 0}(0)=A_{L+1, L+1}(0)$.
We note
$\widehat{I}(\beta)\cup \widehat{I}(-1-\beta)
=\{1,2,\ldots,L+1\}$.

Next, we introduce the two parameters $x$ and $r$ defined as
\begin{eqnarray}
x^{2r}=\frac{q_{2,1}}{q_{1,1}},~~~
r=\left\{\begin{array}{cc}
{\displaystyle \frac{1}{1+\beta}}&
{\rm for}
~~~|\widehat{I}(\beta)|>|\widehat{I}(-1-\beta)|,
\\
{\displaystyle
-\frac{1}{\beta}}&
{\rm for}
~~~|\widehat{I}(\beta)|\leq |\widehat{I}(-1-\beta)|,
\end{array}\right.
\label{parameter:x,r}
\end{eqnarray}
where $|\widehat{I}(\delta)|$ represents the number of 
elements in $\widehat{I}(\delta)$.
By this parametrization, we have
(\ref{thm:Aij(0)}).
From (\ref{assume:p,q1}) and $q_{i,j}\neq 0$, 
we obtain $|x|\neq 0, 1$ and $r>1$.
In this paper, we focus our attention to
\begin{eqnarray}
0<|x|<1,~~~r>1.
\nonumber
\end{eqnarray}
For the case of $|x|>1$, we obtain the same results
under the change $x \to x^{-1}$.

To give the free field construction,
we set $D(k,l; \Phi)$ as
\begin{eqnarray}
&&
D(k,l; {\Phi})=\left\{
\begin{array}{cc}(r-1)\left|\widehat{I}\left(k+1,l+1;-\frac{1}{r},
\widehat{\Phi}\right)\right|
+\left|\widehat{I}\left(k+1,l+1;\frac{1-r}{r},
\widehat{\Phi}\right)\right|&
(0 \leq k \leq l \leq L),\\
0& (0 \leq l<k \leq L),
\end{array}
\right.\nonumber\\
&&
\widehat{I}(k,l;\delta, \widehat{\Phi})
=\{1\leq j \leq L+1|k \leq j \leq l, A_{j-1,j}(0)=\delta\}
~~~(1\leq k \leq l \leq L+1).
\label{def:D(k,l)}
\end{eqnarray}
$D(k,l; \Phi)$ is given by using the matrix 
$(A_{i,j}(0))_{i,j=0}^L$.
The matrix $(A_{i,j}(0))_{i,j=0}^L$ can be constructed from
the Dynkin-diagrams $\Phi$ and $\widehat{\Phi}$,
which we will introduce below.

~\\
{\bf Example}
~~We fix integers $M, N$ $(M \geq N \geq 0, M+N \geq 1)$.
We set $K=1$, $L=M+N+1$, and $j_1=M+1$.
We have
\begin{eqnarray}
&A_{i,i}(0)=\left\{\begin{array}{cc}
\frac{2(r-1)}{r}&{\rm if}~~~1\leq i \leq M,\\
1&{\rm if}~~~i=0, M+1,\\
\frac{~2~}{r}&{\rm if}~~~M+2\leq i \leq L,
\end{array}
\right.\nonumber
\\
&A_{i,i-1}(0)=A_{i-1,i}(0)=\left\{\begin{array}{cc}
\frac{1-r}{r}&{\rm if}~~~1\leq i \leq M+1,\\
-\frac{~1~}{r}&{\rm if}~~~M+2\leq i \leq L+1,
\end{array}
\right.\nonumber
\\
&
A_{k,l}(0)=0~~\left(
|k-l|\geq 2, 1\leq k, l \leq L~~
{\rm or}~~k=0,~l \neq 0, 1, L\right).
\nonumber
\end{eqnarray}
\begin{eqnarray}
\widehat{I}\left(\frac{1-r}{r}\right)=
\{1,2,\ldots,M+1\},
~~~
\widehat{I}\left(\frac{1}{r}\right)=
\{M+2, \ldots,L+1\}.
\nonumber
\end{eqnarray}
We picture $L \times L$ matrix $(A_{i,j}(0))_{i,j=1}^L$
as the standard Dynkin-diagram $\Phi^{st}$ of $A(M,N)$
in Section \ref{Section2}.
We picture $(L+1)\times (L+1)$ matrix $(A_{i,j}(0))_{i,j=0}^L$
as the Dynkin-diagram $\widehat{\Phi}^{st}$ as follows.\\
\begin{center}
\input{Dynkin13.tex}
\\~\\
\end{center}
Here a circle represents an even simple root 
$(\alpha_i, \alpha_i)=2$
and
a crossed circle represents an odd isotropic simple root
$(\alpha_i, \alpha_i)=0$.
The inner product $(\alpha_i,\alpha_j)$ of the roots
and the parameters $A_{i,j}(0)$ correspond
as
$(\alpha_i, \alpha_i)=2 \Leftrightarrow A_{i,i}(0)\neq 1$,
$(\alpha_i, \alpha_i)=0 \Leftrightarrow A_{i,i}(0)=1$,
$(\alpha_i, \alpha_j)=-1 \Leftrightarrow A_{i,j}(0)\neq 0$ $(i\neq j)$.
As additional information, the values of the parameters $A_{j-1,j}(0)$ are written beside the line segment
connecting $\alpha_{j-1}$ and $\alpha_j$.
We have
\begin{eqnarray}
D(0,L; \Phi^{st})=
(N+1)r+M-N.
\nonumber
\end{eqnarray}

~\\
{\bf Example}~~For $L=3$, $K=2$, $j_1=1, j_2=3$, we have
\begin{eqnarray}
\left(A_{i,j}(0)\right)_{i,j=0}^3
=\left(\begin{array}{cccc}
\frac{2(r-1)}{r}&\frac{1-r}{r}&0&\frac{1-r}{r}
\\
\frac{1-r}{r}&1&-\frac{1}{r}&0
\\
0&-\frac{1}{r}&\frac{2}{r}&-\frac{1}{r}
\\
\frac{1-r}{r}&0&-\frac{1}{r}&1
\end{array}\right),~~~
\widehat{I}\left(-\frac{1}{r}\right)=\{2,3\},~~~
\widehat{I}\left(\frac{1-r}{r}\right)=\{1,4\}=\{0,1\}.
\nonumber
\end{eqnarray}
Here we understand subscripts of $A_{i,j}(0)$ with mod.$4$, 
i.e. $A_{0,3}(0)=A_{4,3}(0)$.
We picture $3\times 3$ matrix $(A_{i,j}(0))_{i,j=1}^3$
as nonstandard Dynkin-diagram of $A(1,1)$.
We picture $4\times 4$ matrix $(A_{i,j}(0))_{i,j=0}^3$ as 
the Dynkin-diagram $\widehat{\Phi}$
as follows.
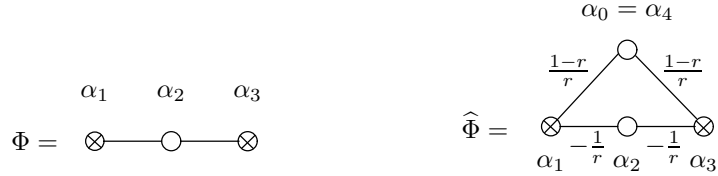
\begin{center}
\input{Dynkin14.tex}~~~~~~~~~~~~~~~\input{Dynkin15.tex}
\end{center}
We have
\begin{eqnarray}
D(0,3;\Phi)=2r,~~~D(1,1;\Phi)=r-1,~~~
D(1,2;\Phi)=2r-2.\nonumber
\end{eqnarray}

\begin{thm}\label{thm:3-1}
Assume (\ref{assume:mutual-locality}), (\ref{assume:q-difference}), 
(\ref{assume:varphiSS}),  (\ref{assume:p,q0}) and (\ref{assume:p,q1}).
Then, up to the freedom \eqref{def:arrange-suffix}, \eqref{def:shift} and \eqref{def:scaling-boson}, 
the parameters
$p_{i,j}, q_{i,j}$, $A_{i,j}(m)$, $s_i(m)$, $\lambda_{i,j}(m)$, $g_i$, and the current $B_j(m)$
are uniquely determined as follows.
Conversely, by choosing these parameters,
(\ref{assume:mutual-locality}), (\ref{assume:q-difference}), and
(\ref{assume:varphiSS}) are satisfied.
\begin{eqnarray}
&&
q_{j,j}=x^{D(1,j-1;\Phi)},~q_{j+1,j}=x^{2r+D(1,j-1;\Phi)}~(1\leq j \leq L),
\nonumber
\\
&&
p_{1,1}=\left\{\begin{array}{cc}
x^2&{\rm if}~~1 \in \widehat{I}\left(-\frac{1}{r}\right),\\
x^{2r-2}&{\rm if}~~1 \in \widehat{I}\left(\frac{1-r}{r}\right),
\end{array}
\right.\nonumber
\\
&&
p_{j,j}=x^{D(1,j-2;\Phi)}\times \left\{
\begin{array}{cc}
x^{r+1}&{\rm if}~~j \in \widehat{I}(-\frac{1}{r}),\\
x^{2r-1}&{\rm if}~~j \in \widehat{I}(\frac{1-r}{r})
\end{array}\right.~(2 \leq j \leq L),
\nonumber\\
&&
p_{j, j-1}=x^{D(1,j-2;\Phi)}\times
\left\{\begin{array}{cc}
x^{2r-2}&{\rm if}~~j \in \widehat{I}(-\frac{1}{r}),\\
x^2&{\rm if}~~j \in \widehat{I}(\frac{1-r}{r})
\end{array}
\right.~(2\leq j \leq L+1),
\nonumber
\\
&&
p_{k,l}=q_{k,l}~~~(k \neq l, l+1, 1\leq k \leq L+1, l \leq l \leq L).
\label{thm:pij}
\end{eqnarray}
\begin{eqnarray}
\begin{array}{c}
s_j(m)=1~~~(m>0, 1\leq j \leq L),\\
s_j(-m)=\left\{\begin{array}{cc}
-1&{\rm if}~~j \in \widehat{J},\\
-\displaystyle{\frac{[m]_x[2(r-1)m]_x}{[rm]_x[(r-1)m]_x}}&
{\rm if}~~j \notin \widehat{J}, j \in I(-\frac{1}{r}),\\
-\displaystyle{
\frac{[(r-1)m]_x[2m]_x}{[rm]_x[m]_x}}&
{\rm if}~~j \notin \widehat{J}, j \in \widehat{I}(\frac{1-r}{r})
\end{array}
\right.
\end{array}~(m>0, 1\leq j \leq L).
\label{thm:sj(m)}
\end{eqnarray}
\begin{eqnarray}
&&
A_{i,i}(0)=\left\{\begin{array}{cc}
1&{\rm if}~j \in \widehat{J},
\\
\displaystyle{\frac{~2~}{r}}&
{\rm if}~~i \notin \widehat{J},~ i \in \widehat{I}(-\frac{1}{r}),
\\
\displaystyle{\frac{2(r-1)}{r}}&
{\rm if}~~i \notin \widehat{J},~ i \in \widehat{I}(\frac{1-r}{r})
\end{array}
\right.~~~(1\leq i \leq L+1),
\nonumber
\\
&&
A_{j-1,j}(0)=A_{j, j-1}(0)=
\left\{\begin{array}{cc}
\displaystyle{-\frac{~1~}{r}}&{\rm if}~~j \in \widehat{I}(-\frac{1}{r}),\\
\displaystyle{\frac{1-r}{r}}&{\rm if}~~j \in \widehat{I}(\frac{1-r}{r})
\end{array}
\right.~(1 \leq j \leq L+1),
\nonumber
\\
&&
A_{k,l}(0)=A_{l,k}(0)=0~(|k-l|\geq 2, 1\leq k<l \leq L~~{\rm or}~~k=0, l \neq 0, 1, L).
\label{thm:Aij(0)}
\end{eqnarray}
\begin{eqnarray}
&&
A_{j,j}(m)=1~(m \neq 0, 1\leq j \leq L),
\nonumber
\\
&&
A_{k,l}(m)=0~(m \neq 0, |k-l|\geq 2, 1\leq k,l \leq L),
\nonumber
\end{eqnarray}
\begin{eqnarray}
&&
\begin{array}{cc}
A_{j-1, j}(m)=\displaystyle{\frac{[m]_x}{[rm]_x}}\times\left\{
\begin{array}{cc}
\displaystyle{\frac{1}{s_{j}(-m)}}&(m>0),\\
\displaystyle{\frac{1}{s_{j-1}(m)}}&(m<0)
\end{array}
\right.&{\rm if}~~j \in \widehat{I}\left(-\frac{1}{r}\right),
\nonumber
\\
A_{j-1,j}(m)=\displaystyle{\frac{[(r-1)m]_x}{[rm]_x}}
\times\left\{
\begin{array}{cc}
\displaystyle{\frac{1}{s_{j}(-m)}}&(m>0),\\
\displaystyle{\frac{1}{s_{j-1}(m)}}&(m<0)
\end{array}
\right.
&{\rm if}~~j \in \widehat{I}\left(\frac{1-r}{r}\right)
\end{array}(2\leq j \leq L),
\nonumber
\\
&&
\begin{array}{c}
A_{j-1,j}(-m)=A_{j,j-1}(m),~~
A_{j,j-1}(-m)=A_{j-1,j}(m)
\end{array}(m>0, 2\leq j \leq L).
\label{thm:Aij(m)}
\end{eqnarray}
\begin{eqnarray}
\lambda_{i,j}(0)=\frac{2r \log x}{D(0,L;\Phi)}
\times
\left\{\begin{array}{cc}
D(0,j-1;\Phi)&{\rm if}~~1\leq j \leq i-1,\\
-D(j,L;\Phi)&
{\rm if}~~i \leq j \leq L
\end{array}
\right.~(1\leq i \leq L+1).
\label{thm:lambda(0)}
\end{eqnarray}
\begin{eqnarray}
\frac{\lambda_{i,j}(m)}{s_j(m)}=\frac{[rm]_x (x-x^{-1})}{
[D(0,L;\Phi)m]_x}
\times
\left\{\begin{array}{cc}
\displaystyle{-x^{(r+D(1,L;\Phi))m}[D(0,j-1;\Phi)m]_x}
&{\rm if}~~1\leq j \leq i-1,\\
\displaystyle{x^{(r-D(0,0;\Phi))m}[D(j,L;\Phi)m]_x}
&
{\rm if}~~i \leq j \leq L
\end{array}
\right.\nonumber\\
(m \neq 0, 1\leq i \leq L+1).
\label{thm:lambda(m)}
\end{eqnarray}
\begin{eqnarray}
g_i=g \times \left\{
\begin{array}{cc}
[r-1]_x
&{\rm if}~~i \in \widehat{I}(-\frac{1}{r}),\\
1&{\rm if}~~i \in \widehat{I}(\frac{1-r}{r})
\end{array}
\right.~(1\leq i \leq L+1).
\label{thm:g}
\end{eqnarray}
\begin{eqnarray}
B_j(z)=g_j\left(\frac{q_{j,j}}{p_{j,j}}-1\right)
:\Lambda_j(z)S_j(q_{j,j}^{-1}z):~(1\leq j \leq L).
\label{thm:B(z)}
\end{eqnarray}
\end{thm}

\begin{prop}
\label{prop:3-2}
The $\Lambda_i(z)$'s satisfy the commutation relations
\begin{eqnarray}
\Lambda_k(z_1)\Lambda_l(z_2)=
\frac{\Theta_{x^{2a}}\left(x^2\frac{z_2}{z_1},~~x^{-2r}\frac{z_2}{z_1},~~x^{2r-2}\frac{z_2}{z_1}\right)}{
\Theta_{x^{2a}}\left(x^{-2}\frac{z_2}{z_1},~~x^{2r}
\frac{z_2}{z_1},~~x^{-2r+2}\frac{z_2}{z_1}\right)
}\Lambda_l(z_2)\Lambda_k(z_1)~~~(1\leq k,l \leq L+1)\,,
\label{eqn:vertex}
\end{eqnarray}
where $a=D(0,L;\Phi)$.
We understand (\ref{eqn:vertex}) in the sense of analytic continuation.
\end{prop}

\begin{prop}
\label{prop:3-3}
The $S_j(w)$'s satisfy the commutation relations
\begin{eqnarray}
&&
S_j(w_1)S_j(w_2)=
S_j(w_2)S_j(w_1)\times
\left\{
\begin{array}{cc}
-1&{\rm if}~~~j \in \widehat{J},\\
\displaystyle{-\left(\frac{w_1}{w_2}\right)^{\frac{2}{r}-1}
\frac{\Theta_{x^{2r}}\left(x^{2}\frac{w_1}{w_2}\right)}{
\Theta_{x^{2r}}\left(x^{2}\frac{w_2}{w_1}\right)}}
&
{\rm if}~~j \notin \widehat{J},~ j \in \widehat{I}(-\frac{1}{r}),
\\
\displaystyle{
-\left(\frac{w_1}{w_2}\right)^{1-\frac{2}{r}}
\frac{\Theta_{x^{2r}}\left(x^{2}\frac{w_2}{w_1}\right)}{
\Theta_{x^{2r}}\left(x^{2}\frac{w_1}{w_2}\right)}}
&
{\rm if}~~j \notin \widehat{J},~ j \in \widehat{I}(\frac{1-r}{r})
\end{array}\right.~~~(1\leq j \leq L),
\nonumber\\
&&
S_{j-1}(w_1)S_j(w_2)=S_j(w_2)S_{j-1}(w_1)
\times
\left\{
\begin{array}{cc}
\displaystyle{
\left(\frac{w_1}{w_2}\right)^{-\frac{1}{r}}
\frac{\Theta_{x^{2r}}\left(x^{r+1}\frac{w_2}{w_1}\right)}{
\Theta_{x^{2r}}\left(x^{r+1}\frac{w_1}{w_2}\right)}}
&
{\rm if}~~j \in \widehat{I}(-\frac{1}{r}),
\\
\displaystyle{
\left(\frac{w_1}{w_2}\right)^{\frac{1}{r}-1}
\frac{\Theta_{x^{2r}}\left(x^{2r-1}\frac{w_2}{w_1}\right)}{
\Theta_{x^{2r}}\left(x^{2r-1}\frac{w_1}{w_2}\right)}}
&
{\rm if}~~j \in \widehat{I}(\frac{1-r}{r})
\end{array}
\right.~~~(2 \leq j \leq L),
\nonumber\\
&&
S_k(w_1)S_l(w_2)=S_l(w_2)S_k(w_1)~~~
(|k-l|\geq 2, 1\leq k, l \leq L).
\label{eqn:screening}
\end{eqnarray}
We understand (\ref{eqn:screening})
in the sense of the analytic continuation.
\end{prop}
In fact, the stronger relation
\begin{eqnarray}
S_j(w_1)S_j(w_2)=(w_1-w_2):S_j(w_1)S_j(w_2):~~~(j \in \widehat{J})
\nonumber
\end{eqnarray}
holds.
This means that the screening currents 
$S_j(w)$ $(j \in \widehat{J})$
are ordinary fermions.

\subsection{Proof of Theorem \ref{thm:3-1}}

In this section, we show Theorem \ref{thm:3-1} and
Proposition \ref{prop:3-3}.
\begin{lem}
\label{lem:3-4}~For $\Lambda_i(z)$ and $S_j(w)$, we obtain
\begin{eqnarray}
&&
\varphi_{\Lambda_i, S_j}(z, w)
=e^{\sum_{k=1}^L \lambda_{i,k}(0)A_{k,j}(0)}
\exp\left(\sum_{k=1}^L \sum_{m=1}^\infty
\frac{1}{m}\lambda_{i,k}(m)A_{k,j}(m)s_j(-m)\left(\frac{w}{z}\right)^m
\right),
\label{normal1}
\\
&&
\varphi_{S_j, \Lambda_i}(w, z)
=\exp\left(\sum_{k=1}^L \sum_{m=1}^\infty
\frac{1}{m}s_j(m)A_{j,k}(m)\lambda_{i,k}(-m)\left(\frac{z}{w}\right)^m
\right)~~~(1\leq i \leq L+1, 1\leq j \leq L),
\label{normal2}
\\
&&
\varphi_{\widetilde{S}_k,\widetilde{S}_l}(w_1, w_2)=
\exp\left(\sum_{m=1}^\infty \frac{1}{m}s_k(m)A_{k,l}(m)s_l(-m)\left(\frac{w_2}{w_1}\right)^m \right)~~~(1\leq k, l \leq L).
\label{normal3}
\end{eqnarray}

Assume (\ref{assume:mutual-locality}), we obtain
\begin{eqnarray}
\varphi_{\Lambda_k, \Lambda_l}(z_1, z_2)=
\exp\left(\sum_{i=1}^L\sum_{m=1}^\infty
\frac{1}{m}\frac{\lambda_{k,i}(m)}{s_i(m)}
(q_{l,i}^{-m}-p_{l,i}^{-m})
\left(\frac{z_2}{z_1}\right)^m \right)~~~(1\leq k, l \leq L+1).
\label{normal4}
\end{eqnarray}
\end{lem}

\begin{lem}
\label{lem:3-5}~
Mutual locality (\ref{assume:mutual-locality}) holds,
if and only if (\ref{eqn:lambda(0)})
and (\ref{eqn:lambda(m)}) are satisfied.
\begin{eqnarray}
&&
\sum_{k=1}^L
\lambda_{i,k}(0)A_{k,j}(0)=\log\left(\frac{q_{i,j}}{p_{i,j}}\right)
~~~(1\leq i \leq L+1, 1\leq j \leq L),
\label{eqn:lambda(0)}
\\
&&
\sum_{k=1}^L
\lambda_{i,k}(m)A_{k,j}(m)s_j(-m)=q_{i,j}^m-p_{i,j}^m
~~~(m \neq 0, 1\leq i \leq L+1, 1\leq j \leq L).
\label{eqn:lambda(m)}
\end{eqnarray}
\end{lem}
{\it Proof of Lemmas \ref{lem:3-4} and \ref{lem:3-5}.}~
Using the standard formula
\begin{eqnarray}
e^Ae^B=e^{[A,B]}e^Be^A~~~([[A,B],A]=0~{\rm and}~[[A,B],B]=0),
\nonumber
\end{eqnarray}
we obtain (\ref{normal1}), (\ref{normal2}), (\ref{normal3}), and 
\begin{eqnarray}
\varphi_{\Lambda_k, \Lambda_l}(z_1, z_2)=
\exp\left(\sum_{i, j=1}^L \sum_{m=1}^\infty
\frac{1}{m}\lambda_{k,i}(m)A_{i,j}(m)\lambda_{l,j}(-m)\left(\frac{z_2}{z_1}\right)^m \right)~~~(1\leq k, l \leq L+1).
\label{normal5}
\end{eqnarray}

Considering (\ref{normal1}), (\ref{normal2}), and the expansions
\begin{eqnarray}
&&
{\displaystyle \frac{w-p_{i,j}^{-1} z}{w-q_{i,j}^{-1} z}}=
\exp\left(\log\left(\frac{q_{i,j}}{p_{i,j}}\right)-\sum_{m=1}^\infty \frac{1}{m}(p_{i,j}^m-q_{i,j}^m)\left(\frac{w}{z}\right)^m\right)
~~~(|z| \gg |w|),
\label{expansion1}\\
&&
\frac{w-p_{i,j}^{-1} z}{w-q_{i,j}^{-1} z}=
\exp\left(-\sum_{m=1}^\infty \frac{1}{m}(p_{i,j}^{-m}-q_{i,j}^{-m})\left(\frac{z}{w}\right)^{m}\right)
~~~(|w| \gg |z|),
\label{expansion2}
\end{eqnarray}
we obtain
(\ref{eqn:lambda(0)}) and (\ref{eqn:lambda(m)})
from (\ref{assume:mutual-locality}).
Substituting (\ref{eqn:lambda(m)}) for (\ref{normal5}),
we have (\ref{normal4}).

Conversely, if we assume (\ref{eqn:lambda(0)}) and 
(\ref{eqn:lambda(m)}), we obtain (\ref{assume:mutual-locality})
from (\ref{normal1}), (\ref{normal2}), (\ref{expansion1}), and (\ref{expansion2}).
~~$\Box$

From the linear equations (\ref{eqn:lambda(0)}) and (\ref{eqn:lambda(m)}),
$\lambda_{i,j}(m)$ are expressed in terms of the other parameters.

\begin{lem}\label{lem:3-6}~
We assume (\ref{assume:mutual-locality}) and (\ref{assume:p,q0}).
The commutativity (\ref{assume:q-difference}) holds,
if and only if (\ref{assume:start0}), (\ref{assume:start1}), (\ref{assume:start2}), and (\ref{def:current-B}) are satisfied.
\begin{eqnarray}
&&
p_{k,l}=q_{k,l}~~~(k\neq l, l+1, 1\leq k \leq L+1, 1\leq l \leq L),
\label{assume:start0}
\\
&&
q_{k,k}^{\frac{1}{2}A_{k,k}(0)}:\Lambda_k(z)
S_k
\left(q_{k,k}^{-1}z\right):=
q_{k+1,k}^{\frac{1}{2}A_{k,k}(0)}:\Lambda_{k+1}(z)
S_k\left(q_{k+1,k}^{-1} z \right):
~~~(1\leq k \leq L),
\label{assume:start1}
\\
&&
\frac{g_{k+1}}{g_k}=-\left(\frac{q_{k+1,k}}{q_{k,k}}\right)^{\frac{1}{2}A_{k,k}(0)}
\frac{\frac{q_{k,k}}{p_{k,k}}-1}{
\frac{~q_{k+1,k}}{p_{k+1,k}}-1~}~~~(1\leq k \leq L),
\label{assume:start2}
\\
&&
B_k(z)=g_k\left(\frac{q_{k,k}}{p_{k,k}}-1\right)
:\Lambda_k(z)S_k(q_{k,k}^{-1}z):~~~(1\leq k \leq L).
\label{def:current-B}
\end{eqnarray}
\end{lem}
{\it Proof of Lemma \ref{lem:3-6}.}~
From (\ref{assume:mutual-locality}), we obtain
\begin{eqnarray}
[\Lambda_i(z),S_j(w)]=\left(\frac{q_{i,j}}{p_{i,j}}-1\right)
\delta\left(\frac{q_{i,j}w}{z}\right):\Lambda_i(z)S_j(q_{i,j}^{-1}z):
~~(1\leq i \leq L+1, 1\leq j \leq L).\label{eqn:lambda-S}
\end{eqnarray}

Considering (\ref{assume:p,q0}) and (\ref{eqn:lambda-S}), we know that
(\ref{assume:q-difference}) holds, 
if and only if (\ref{assume:start0}) and
\begin{eqnarray}
B_j(z)=g_j\left(\frac{q_{j,j}}{p_{j,j}}-1\right):\Lambda_j(z)S_j(q_{j,j}^{-1}z):
=-g_{j+1}\left(\frac{q_{j+1,j}}{p_{j+1,j}}-1\right):\Lambda_{j+1}(z)S_j(q_{j+1,j}^{-1}z):~(1\leq j \leq L)
\label{eqn:current-B}
\end{eqnarray}
are satisfied.
(\ref{eqn:current-B}) holds,
if and only if (\ref{assume:start1}), (\ref{assume:start2}), and (\ref{def:current-B}) are satisfied.
Hence, we obtain this lemma.~~~
$\Box$

We use the abbreviation $h_{k,l}(w)$ $(1\leq k,l \leq L)$ as
\begin{eqnarray}
h_{k,l}\left(\frac{w_2}{w_1}\right)=
\varphi_{\widetilde{S}_k, \widetilde{S}_l}(w_1, w_2).
\label{def:h}
\end{eqnarray}

\begin{lem}\label{lem:3-7}~
We assume (\ref{assume:mutual-locality}) and (\ref{assume:start1}).
Then, $h_{k,l}(w)$ in (\ref{def:h}) satisfy the $q$-difference equations
\begin{eqnarray}
&&
\begin{array}{c}
{\displaystyle
\frac{w-p_{k,k}^{-1}}{w-q_{k,k}^{-1}}
h_{k,k}\left(q_{k,k}^{-1}w\right)=
\frac{w-p_{k+1,k}^{-1}}{w-q_{k+1, k}^{-1}}
h_{k,k}\left(q_{k+1,k}^{-1}w \right)},
\\
{\displaystyle
\left(\frac{q_{k+1,k}}{q_{k,k}}\right)^{A_{k,k}(0)-1}
\frac{p_{k+1,k}}{p_{k,k}}
\frac{1-p_{k,k}w}{1-q_{k,k}w}h_{k,k}\left(q_{k,k}w\right)=
\frac{1-p_{k+1,k}w}{1-q_{k+1,k}w}h_{k,k}\left(q_{k+1,k}w\right)}
\end{array}
~(1\leq k \leq L),
\label{eqn:hii-1}
\end{eqnarray}
and
\begin{eqnarray}
\begin{array}{c}
{\displaystyle \frac{h_{k,k+1}(q_{k,k}w)}{h_{k,k+1}(q_{k+1,k}w)}
=\frac{q_{k+1,k+1}}{p_{k+1,k+1}}
\left(\frac{q_{k,k}}{q_{k+1,k}}\right)^{A_{k,k+1}(0)}
\frac{1-p_{k+1,k+1}w}{1-q_{k+1,k+1}w}},
\\
{\displaystyle 
\frac{h_{k,k+1}\left(q_{k+1,k+1}^{-1}w \right)}
{h_{k,k+1}\left(q_{k+2,k+1}^{-1} w \right)}
=\frac{1-q_{k+1,k}^{-1} w}{1-p_{k+1,k}^{-1}w}},
\\
{\displaystyle
\frac{h_{k+1,k}(q_{k+2,k+1}w)}{h_{k+1,k}(q_{k+1,k+1}w)}
=\frac{q_{k+1,k}}{p_{k+1,k}}
\left(\frac{q_{k+2,k+1}}{q_{k+1,k+1}}\right)^{A_{k,k+1}(0)}
\frac{1-p_{k+1,k}w}{1-q_{k+1,k}w}},
\\
{\displaystyle
\frac{h_{k+1,k}\left(q_{k,k}^{-1} w \right)}{h_{k+1,k}
\left(q_{k+1,k}^{-1} w \right)}
=\frac{1-p_{k+1,k+1}^{-1}w}{1-q_{k+1,k+1}^{-1} w}}
\end{array}
~(1\leq k \leq L-1).
\label{eqn:q-diff1}
\end{eqnarray}
\end{lem}
{\it Proof of Lemma \ref{lem:3-7}.}~
Multiplying (\ref{assume:start1}) by the screening currents on the left or right and considering the normal orderings, we obtain
(\ref{eqn:hii-1}) and (\ref{eqn:q-diff1})
as necessary conditions.~~~
$\Box$

\begin{lem}\label{lem:3-8}~
The relations (\ref{eqn:p,q0}) and (\ref{eqn:p,q1}) hold, 
if (\ref{assume:mutual-locality}), (\ref{assume:varphiSS}), (\ref{assume:p,q1}), and (\ref{assume:start1}) are satisfied.
\begin{eqnarray}
\begin{array}{c}
q_{k+1,k+1}=q_{k,k}x^{(1+A_{k,k+1}(0))r},\\
q_{k+1,k}=q_{k,k}x^{2r},
\\
p_{k+1,k+1}=q_{k,k}x^{(1-A_{k,k+1}(0))r},
\\
p_{k+1,k}=q_{k,k}x^{2(1+A_{k,k+1}(0))r}
\end{array}~~~(1\leq k \leq L-1).
\label{eqn:p,q0}
\\
\begin{array}{cc}
\begin{array}{c}
p_{k,k}=q_{k,k}x^{A_{k,k}(0)r},
\\
p_{k+1,k}=q_{k,k}x^{(2-A_{k,k}(0))r}
\end{array}
&
{\rm if}~~~k \notin \widehat{J}, 
\\
p_{k+1,k}=p_{k,k}
&
{\rm if}~~~k \in \widehat{J}
\end{array}
~~~
(1\leq k \leq L).
\label{eqn:p,q1}
\end{eqnarray}
\end{lem}

\begin{lem}
\label{lem:3-9}
The relation (\ref{eqn:A(m)}) holds, 
if (\ref{assume:mutual-locality}), (\ref{assume:varphiSS}), (\ref{assume:p,q1}), and (\ref{assume:start1}) are satisfied.
\begin{eqnarray}
&&
s_k(m)s_k(-m)=
\left\{
\begin{array}{cc}
-1&{\rm if}~~~k \in \widehat{J},
\\
\displaystyle{
-\frac{[\frac{1}{2}A_{k,k}(0)rm]_x[(2-A_{k,k}(0))r m]_x}{
[\frac{1}{2}(2-A_{k,k}(0))rm]_x[rm]_x}}
&{\rm if}~~~k \notin \widehat{J}
\end{array}\right.~(m>0, 1\leq k \leq L),
\nonumber
\\
&&
\begin{array}{c}
s_k(m)A_{k,k+1}(m)s_{k+1}(-m)
=\displaystyle{-\frac{[A_{k,k+1}(0) r m]_x}{[rm]_x}},
\\
s_{k+1}(m)A_{k+1,k}(m)s_{k}(-m)
=\displaystyle{-\frac{[A_{k+1,k}(0) r m]_x}{[rm]_x}}
\end{array}
~~(m>0, 1\leq k \leq L-1),
\nonumber\\
&&
A_{k,l}(m)=0~~~(m>0, |k-l|\geq 2, 1\leq k,l \leq L).
\label{eqn:A(m)}
\end{eqnarray}
\end{lem}
{\it Proof of Lemmas \ref{lem:3-8} and \ref{lem:3-9}.}~
From Lemma \ref{lem:3-7}, we obtain
the $q$-difference equations 
(\ref{eqn:hii-1}) and (\ref{eqn:q-diff1}).
From (\ref{normal3}) and (\ref{def:h}), 
the constant term of $h_{k,l}(w)$ is 1.
Comparing the Taylor expansions for both sides of 
(\ref{eqn:hii-1}) and (\ref{eqn:q-diff1}), 
we obtain
\begin{eqnarray}
&&
\frac{p_{k+1,k}}{p_{k,k}}
\left(\frac{q_{k+1,k}}{q_{k,k}}\right)^{A_{k,k}(0)-1}
=1~~~(1\leq k \leq L),
\label{eqn:p,q6}
\\
&&
\frac{q_{k+1,k+1}}{p_{k+1,k+1}}
\left(\frac{q_{k,k}}{q_{k+1,k}}\right)^{A_{k,k+1}(0)}=1,~~
\frac{q_{k+1,k}}{p_{k+1,k}}
\left(\frac{q_{k+2,k+1}}{q_{k+1,k+1}}\right)^{A_{k,k+1}(0)}
=1~~~(1\leq k \leq L-1).\label{eqn:p,q2}
\end{eqnarray}

First, we study the $q$-difference equations
in (\ref{eqn:q-diff1}).
Upon the specialization (\ref{eqn:p,q2}), 
we obtain solutions of
(\ref{eqn:q-diff1}) as
\begin{eqnarray}
h_{k,k+1}(w)=\exp\left(-\sum_{m=1}^\infty
\frac{1}{m}\frac{\left(\frac{p_{k+1,k+1}}{q_{k,k}}\right)^m-
\left(\frac{q_{k+1,k+1}}{q_{k,k}}\right)^m}{1-\left(\frac{q_{k+1,k}}{q_{k,k}}\right)^m}w^m\right)
=
\exp\left(-\sum_{m=1}^\infty
\frac{1}{m}\frac{\left(\frac{q_{k+2,k+1}}{p_{k+1,k}}\right)^m-
\left(\frac{q_{k+2,k+1}}{q_{k+1,k}}\right)^m}{1-\left(\frac{q_{k+2,k+1}}{q_{k+1,k+1}}\right)^m}w^m\right),
\label{eqn:h12-1}
\\
h_{k+1,k}(w)=
\exp\left(-\sum_{m=1}^\infty
\frac{1}{m}\frac{\left(\frac{q_{k+1,k}}{q_{k+1,k+1}}\right)^m-
\left(\frac{p_{k+1,k}}{q_{k+1,k+1}}\right)^m}{1-
\left(\frac{q_{k+2,k+1}}{q_{k+1,k+1}}\right)^m}w^m\right)
=
\exp\left(-\sum_{m=1}^\infty
\frac{1}{m}
\frac{\left(\frac{q_{k+1,k}}{q_{k+1,k+1}}\right)^m-
\left(\frac{q_{k+1,k}}{p_{k+1,k+1}}\right)^m}
{1-\left(\frac{q_{k+1,k}}{q_{k,k}}\right)^m}
w^m\right).
\label{eqn:h21-1}
\end{eqnarray}
Here we used $|q_{k+1,k}/q_{k,k}|\neq 1$ $(1\leq k \leq L)$ 
assumed in (\ref{assume:p,q1}).
From the compatibility of the two formulae for $h_{k,k+1}(w)$
in (\ref{eqn:h12-1})
[ or $h_{k+1,k}(w)$ in 
(\ref{eqn:h21-1})],
there are two possible choices for 
$q_{k,k}$, $q_{k+1,k+1}$, $q_{k+1,k}$,
and $q_{k+2,k+1}$.
\begin{eqnarray}
\mathrm{(i)}~~\frac{q_{k+1,k}}{q_{k,k}}=\frac{q_{k+2,k+1}}{q_{k+1,k+1}}~~~{\rm or}~~~
\mathrm{(ii)}~~\frac{q_{k+1,k}}{q_{k,k}}=\frac{q_{k+1,k+1}}{q_{k+2,k+1}}.
\label{choice:q}
\end{eqnarray}
 
First, we consider the case of 
$\mathrm{(ii)}$~$
q_{k+1,k}/q_{k,k}=q_{k+1,k+1}/q_{k+2,k+1}$ in (\ref{choice:q}).
From the compatibility of the two formulae for $h_{k,k+1}(w)$
in (\ref{eqn:h12-1})
[and $h_{k+1,k}(w)$ in (\ref{eqn:h21-1})],
we obtain
\begin{eqnarray}
\left(\frac{p_{k+1,k+1}}{q_{k,k}}\right)^m+\left(\frac{q_{k+1,k+1}}{
p_{k+1,k}}\right)^m=
\left(\frac{q_{k+1,k+1}}{q_{k+1,k}}\right)^m+
\left(\frac{q_{k+1,k+1}}{q_{k,k}}\right)^m~~(m \neq 0).
\label{eqn:p,q4}
\end{eqnarray}
From (\ref{eqn:p,q4}) for $m=1,2$, we obtain 
$p_{k+1,k+1}/p_{k+1,k}=q_{k+1,k+1}/q_{k+1,k}$.
Combining (\ref{eqn:p,q4}) for $m=1$ 
and $p_{k+1,k+1}/p_{k+1,k}=q_{k+1,k+1}/q_{k+1,k}$, we obtain
$q_{k,k}=p_{k+1,k}$ or $q_{k+1,k+1}=p_{k+1,k+1}$.
For the case of $q_{k,k}=p_{k+1,k}$, 
we obtain $A_{k,k+1}(0)=1$ from (\ref{eqn:p,q2}).
For the case of $q_{k+1,k+1}=p_{k+1,k+1}$, 
we obtain $A_{k,k+1}(0)=0$ from (\ref{eqn:p,q2}).
$A_{k,k+1}(0)=0$ and $A_{k,k+1}(0)=1$ contradict with 
$-1<A_{k,k+1}(0)<0$ assumed in (\ref{assume:p,q1}).
Hence, the case of $\mathrm{(ii)}~q_{k+1,k}/q_{k,k}=q_{k+1,k+1}/
q_{k+2,k+1}$ is impossible.

Next, we consider the case of 
$\mathrm{(i)}$~
$q_{k+1,k}/q_{k,k}=q_{k+2,k+1}/q_{k+1,k+1}$ in (\ref{choice:q}).
From exclusion of the case $\mathrm{(ii)}$ and the parametrization
(\ref{parameter:x,r}), we can parametrize
\begin{eqnarray}
\frac{q_{2,1}}{q_{1,1}}=\frac{q_{3,2}}{q_{2,2}}=\cdots=
\frac{q_{L+1,L}}{q_{L,L}}=x^{2r}.
\label{eqn:p,q8}
\end{eqnarray}
From the compatibility of the two formulae for $h_{k,k+1}(w)$
in (\ref{eqn:h12-1})
[ and $h_{k+1, k}(w)$
in (\ref{eqn:h21-1})],
we obtain
\begin{eqnarray}
&&
h_{k,k+1}(w)=\exp\left(-\sum_{m=1}^\infty
\frac{1}{m}\frac{[A_{k,k+1}(0)rm]_x}{[rm]_x}x^{-(A_{k,k+1}(0)+1)rm}
\left(\frac{q_{k+1,k+1}}{q_{k,k}}\right)^mw^m
\right),
\label{eqn:h12-3}
\\
&&
h_{k+1,k}(w)=\exp\left(-\sum_{m=1}^\infty
\frac{1}{m}\frac{[A_{k,k+1}(0)rm]_x}{[rm]_x}x^{(A_{k,k+1}(0)+1)rm}
\left(\frac{q_{k,k}}{q_{k+1,k+1}}\right)^mw^m
\right).
\label{eqn:h21-3}
\end{eqnarray}
We used (\ref{eqn:p,q2}) and (\ref{eqn:p,q8}).
From $h_{k,k+1}(w)=h_{k+1,k}(w)$ assumed in (\ref{assume:varphiSS}), we obtain
\begin{eqnarray}
\frac{q_{k+1,k+1}}{q_{k,k}}=x^{(A_{k,k+1}(0)+1)r}.
\label{eqn:p,q3}
\end{eqnarray}
Considering (\ref{eqn:p,q2}), (\ref{eqn:p,q8}), and (\ref{eqn:p,q3}), 
we obtain (\ref{eqn:p,q0}).
From (\ref{eqn:h12-3}), (\ref{eqn:h21-3}), and
(\ref{eqn:p,q3}), we obtain
\begin{eqnarray}
h_{k,k+1}(w)=h_{k+1,k}(w)=
\exp\left(-\sum_{m=1}^\infty
\frac{1}{m}\frac{[A_{k,k+1}(0)rm]_x}{[rm]_x}w^m
\right).
\label{sym:h}
\end{eqnarray}
Considering (\ref{normal3}), (\ref{def:h}), and (\ref{sym:h}), we obtain the second half of (\ref{eqn:A(m)}).

Next, we study the $q$-difference equations in (\ref{eqn:hii-1}).
Upon the specialization (\ref{eqn:p,q6}), the compatibility condition of the equations in
(\ref{eqn:hii-1}) is
\begin{eqnarray}
(p_{k,k}-p_{k+1,k})(p_{k,k}p_{k+1,k}-q_{k,k}q_{k+1,k})=0
~~~(1\leq k \leq L).
\label{eqn:p,q7}
\end{eqnarray}

First, we study the case of $A_{k,k}(0)=1$.
We obtain $p_{k,k}=p_{k+1,k}$ 
in the second half of (\ref{eqn:p,q1})
from (\ref{eqn:p,q6}).
Solving (\ref{eqn:hii-1}) upon $p_{k,k}=p_{k+1,k}$, 
we obtain $h_{k,k}(w)=1-w$.
Considering (\ref{normal3}) and (\ref{def:h}), we obtain 
$s_k(m)s_k(-m)=-1$ $(m>0)$ in the first half of (\ref{eqn:A(m)}).

Next, we study the case of $A_{k,k}(0)\neq 1$.
We obtain $p_{k,k}\neq p_{k+1,k}$ from (\ref{assume:p,q1}) and (\ref{eqn:p,q6}).
Then, we obtain 
$p_{k,k}p_{k+1,k}=q_{k,k}q_{k+1,k}$ from (\ref{eqn:p,q7}).
Combining $p_{k,k}p_{k+1,k}=q_{k,k}q_{k+1,k}$ 
and (\ref{eqn:p,q6}), we obtain (\ref{eqn:p,q1}).
Solving (\ref{eqn:hii-1}), we obtain
\begin{eqnarray}
h_{k,k}(w)=\exp\left(
-\sum_{m=1}^\infty \frac{1}{m}
{\displaystyle
\frac{\left[\frac{1}{2}A_{k,k}(0)rm\right]_x\left[(2-A_{k,k}(0))rm\right]_x}{
\left[\frac{1}{2}(2-A_{k,k}(0))r m\right]_x [rm]_x}}w^m
\right).\nonumber
\end{eqnarray}
We used 
$|q_{k+1,k}/q_{k,k}|\neq 1$ $(1\leq k \leq L)$ in (\ref{assume:p,q1})
and
$q_{k+1,k}/q_{k,k}=x^{2r}$ $(1\leq k \leq L)$ in (\ref{eqn:p,q8}).
Considering (\ref{normal3}) and (\ref{def:h}), we obtain the first half of (\ref{eqn:A(m)}).
~~$\Box$

\begin{lem}
\label{lem:3-10}
The relation (\ref{eqn:A(0)-1}) holds for $(A_{i,j}(0))_{i,j=0}^L$,
if (\ref{assume:mutual-locality}), (\ref{assume:q-difference}), (\ref{assume:varphiSS}), (\ref{assume:p,q0}), and (\ref{assume:p,q1}) are satisfied.
\end{lem}
{\it Proof of Lemma \ref{lem:3-10}.}~
We obtain $A_{k,l}(0)$ $(|k-l|\geq 2, 1\leq k,l \leq L)$
from (\ref{assume:varphiSS}).
From Lemma \ref{lem:3-6}, we have (\ref{assume:start1}).
From Lemma \ref{lem:3-8}, we have
(\ref{eqn:p,q0}) and (\ref{eqn:p,q1}).
From the compatibility of
(\ref{eqn:p,q0}) and (\ref{eqn:p,q1}), we obtain
the following relations for $(A_{i,j}(0))_{i,j=1}^L$.
\begin{eqnarray}
&&
A_{k+1,k}(0)=\left\{\begin{array}{cc}
A_{k-1,k}(0)&{\rm if}~~k \notin \widehat{J},\\
-1-A_{k-1,k}(0)&{\rm if}~~k \in \widehat{J}
\end{array}
\right.~~~(2\leq k \leq L-1),\nonumber\\
&&
A_{k,k}(0)=\left\{\begin{array}{cc}
-2A_{k-1,k}(0)&{\rm if}~~k \notin \widehat{J},\\
1&{\rm if}~~k \in \widehat{J}
\end{array}
\right.~(2\leq k \leq L),~~~
A_{1,1}(0)=
\left\{\begin{array}{cc}
-2A_{1,2}(0)&{\rm if}~~1 \notin \widehat{J},\\
1&{\rm if}~~1 \in \widehat{J},
\end{array}
\right.
\nonumber\\
&&
A_{k,l}(0)=0~~~
(|k-l|\geq 2, 1\leq k,l \leq L).
\end{eqnarray}
Solving these equations, we obtain (\ref{eqn:A(0)-1})
for $1\leq i \leq L$, $2\leq j \leq L$, and $1\leq k,l \leq L$.
The extension to $(A_{i,j}(0))_{i,j=0}^L$ is direct consequence of
the definition ({\ref{eqn:A(0)-0}}).
~$\Box$

\begin{prop}\label{prop:3-11}~
The relations (\ref{assume:mutual-locality}), (\ref{assume:q-difference}), and (\ref{assume:varphiSS}) hold,
if the parameters $p_{i,j}, q_{i,j}$, $A_{i,j}(m)$, $s_i(m)$,
$g_i$, and $\lambda_{i,j}(m)$ are determined by 
(\ref{thm:Aij(0)}), (\ref{eqn:lambda(0)}), (\ref{eqn:lambda(m)}), (\ref{assume:start0}), (\ref{assume:start2}), 
(\ref{eqn:p,q0}), (\ref{eqn:p,q1}), and (\ref{eqn:A(m)}).
\end{prop}
{\it Proof of Proposition \ref{prop:3-11}.}~
First, we will show the formulae
(\ref{thm:pij}), (\ref{thm:sj(m)}),
(\ref{thm:Aij(m)}), (\ref{thm:lambda(0)}), (\ref{thm:lambda(m)}),
and (\ref{thm:g}) in Theorem \ref{thm:3-1}.
Let $q_{1,1}=s$.
From (\ref{thm:Aij(0)}), (\ref{assume:start0}),
(\ref{eqn:p,q0}), and (\ref{eqn:p,q1}), we have
$q_{j,j}=s x^{D(1,j-1;\Phi)}$,~~
$q_{j+1,j}=s x^{2r+D(1,j-1;\Phi)}~(1\leq j \leq L)$,
~
$p_{1,1}=s \times \left\{\begin{array}{cc}
x^2&{\rm if}~~1 \in \widehat{I}\left(-\frac{1}{r}\right),\\
x^{2r-2}&{\rm if}~~1 \in \widehat{I}\left(\frac{1-r}{r}\right)
\end{array}
\right.$,~
$p_{j,j}=s x^{D(1,j-2;\Phi)}\times \left\{
\begin{array}{cc}
x^{r+1}&{\rm if}~~j \in \widehat{I}(-\frac{1}{r}),\\
x^{2r-1}&{\rm if}~~j \in \widehat{I}(\frac{1-r}{r})
\end{array}\right.$ $(2 \leq j \leq L)$,~
$p_{j, j-1}=s x^{D(1,j-2;\Phi)}\times
\left\{\begin{array}{cc}
x^{2r-2}&{\rm if}~~j \in \widehat{I}(-\frac{1}{r}),\\
x^2&{\rm if}~~j \in \widehat{I}(\frac{1-r}{r})
\end{array}
\right.$ 
$(2\leq j \leq L+1)$,
$p_{k,l}=q_{k,l}$ $(k \neq l, l+1, 1\leq k \leq L+1, l \leq l \leq L)$.
Upon the specialization $s=1$, we have (\ref{thm:pij}).
From (\ref{thm:pij}), (\ref{thm:Aij(0)}), and (\ref{assume:start2}),
we have (\ref{thm:g}).
From (\ref{eqn:A(m)}), we have
\begin{eqnarray}
&&
s_k(-m)=-\frac{1}{\alpha_k(m)}
\times
\left\{\begin{array}{cc}
1&{\rm if}~~k \in \widehat{J},\\
\displaystyle{\frac{[\frac{1}{2}A_{k,k}(0)rm]_x
[(2-A_{k,k}(0))rm]_x}{[\frac{1}{2}(2-A_{k,k}(0))rm]_x[rm]_x}}&{\rm if}~~k \notin \widehat{J}
\end{array}~~(m>0),
\right.\nonumber
\\
&&
A_{k\pm 1,k}(m)=
\frac{\alpha_k(m)}{\alpha_{k\pm 1}(m)}
\frac{[A_{k\pm 1,k}(0)rm]_x}{[rm]_x}
\times
\left\{\begin{array}{cc}
1&{\rm if}~~k \in \widehat{J},
\\
\displaystyle{
\frac{[rm]_x[\frac{1}{2}(2-A_{k,k}(0))rm]_x}{
[\frac{1}{2}A_{k,k}(0)rm]_x[(2-A_{k,k}(0))rm]_x}}
&{\rm if}~~k \notin \widehat{J}
\end{array}~~(m>0),
\right.
\nonumber
\\
&&
A_{k,k\pm1}(-m)=A_{k\pm1,k}(m)~~(m>0).
\nonumber
\end{eqnarray}
Here the signs of the formulae are in the same order.
Here we set $s_k(m)=\alpha_k(m)$ $(m>0, 1\leq k \leq L)$.
Setting $\alpha_k(m)=1$ $(m>0, 1\leq k \leq L)$ provides
(\ref{thm:sj(m)}) and (\ref{thm:Aij(m)}).
Solving (\ref{eqn:lambda(0)}) and (\ref{eqn:lambda(m)}),
we obtain $\lambda_{i,j}(m)$ in (\ref{thm:lambda(0)}) and (\ref{thm:lambda(m)}).
Solving (\ref{eqn:lambda(0)}) and (\ref{eqn:lambda(m)}) for arbitrary $\alpha_k(m)$,
we obtain $\lambda_{i,j}(m)\alpha_j(m)$.
Now we obtained the formulae
(\ref{thm:pij}), (\ref{thm:sj(m)}),
(\ref{thm:Aij(m)}), (\ref{thm:lambda(0)}), (\ref{thm:lambda(m)}),
and (\ref{thm:g}). 
As a by-product of calculation, we proved that
there is no indeterminacy in the free field realization except for
(\ref{def:arrange-suffix}),
(\ref{def:shift}),
and (\ref{def:scaling-boson}), which is part of 
Theorem \ref{thm:3-1}.

Next, we will derive
(\ref{assume:mutual-locality}), 
(\ref{assume:q-difference}), 
and (\ref{assume:varphiSS}).
From (\ref{thm:sj(m)}) and (\ref{thm:Aij(m)}),
we obtain the symmetry (\ref{assume:varphiSS}) by direct calculation.
Because $\lambda_{i,j}(m)$ are determined by
(\ref{eqn:lambda(0)}) and (\ref{eqn:lambda(m)}),
the mutual locality (\ref{assume:mutual-locality})
holds from Lemma \ref{lem:3-5}.
From (\ref{thm:sj(m)}) and (\ref{thm:lambda(m)}),
we have (\ref{assume:start1}) by direct calculation.
From (\ref{assume:mutual-locality}), (\ref{assume:start0}),
(\ref{assume:start1}), and (\ref{assume:start2}),
we obtain 
$[T_1(z),S_j(w)]=\left(\frac{q_{j,j}}{p_{j,j}}-1\right):\Lambda_j(z)S_j(q_{j,j}^{-1}z):
\left(\delta\left(\frac{q_{j,j}w}{z}\right)-
\delta\left(\frac{q_{j+1,j}w}{z}\right)\right)$~$(1\leq j \leq L)$.
Hence, we have the commutativity (\ref{assume:q-difference})
upon the condition (\ref{thm:B(z)}).
We derived
(\ref{assume:mutual-locality}), 
(\ref{assume:q-difference}), 
and (\ref{assume:varphiSS}).~~$\Box$

~\\
{\it Proof of Theorem \ref{thm:3-1}.}~
We assume the relations
(\ref{assume:mutual-locality}), 
(\ref{assume:q-difference}), 
(\ref{assume:varphiSS}), 
(\ref{assume:p,q0}), 
and (\ref{assume:p,q1}).
From Lemmas \ref{lem:3-5}, \ref{lem:3-6}, \ref{lem:3-8},
\ref{lem:3-9}, and \ref{lem:3-10}, we obtain
the relations
(\ref{thm:Aij(0)}), (\ref{eqn:lambda(0)}), (\ref{eqn:lambda(m)}), (\ref{assume:start0}), (\ref{assume:start2}), 
(\ref{eqn:p,q0}), (\ref{eqn:p,q1}), and (\ref{eqn:A(m)}).
In proof of Proposition \ref{prop:3-11},
we have obtained
$p_{i,j}$, $q_{i,j}$, $s_j(m)$, $A_{i,j}(m)$, $g_i$, $\lambda_{i,j}(m)$,
and $B_j(z)$ in
(\ref{thm:pij}), (\ref{thm:sj(m)}), (\ref{thm:Aij(m)}), 
(\ref{thm:lambda(0)}), (\ref{thm:lambda(m)}), 
(\ref{thm:g}), and (\ref{thm:B(z)})
from the relations (\ref{thm:Aij(0)}), (\ref{eqn:lambda(0)}), (\ref{eqn:lambda(m)}), (\ref{assume:start0}), (\ref{assume:start2}), 
(\ref{eqn:p,q0}), (\ref{eqn:p,q1}), and (\ref{eqn:A(m)}).
Moreover,
in proof of Proposition \ref{prop:3-11},
we have proved that
there is no indeterminacy in the free field realization except for
(\ref{def:arrange-suffix}),
(\ref{def:shift}),
and (\ref{def:scaling-boson}).

Conversely,
in proof of Proposition \ref{prop:3-11}, we have proved that
the relations 
(\ref{assume:mutual-locality}), (\ref{assume:q-difference}), 
and (\ref{assume:varphiSS}) hold, if
the relations (\ref{thm:pij}), (\ref{thm:sj(m)}), 
(\ref{thm:Aij(0)}), (\ref{thm:Aij(m)}), 
(\ref{thm:lambda(0)}), (\ref{thm:lambda(m)}), 
(\ref{thm:g}), and (\ref{thm:B(z)}) are satisfied.~~$\Box$

~\\
{\it Proof of Proposition \ref{prop:3-3}.}~
Using $h_{k,l}(w)$ in (\ref{def:h}), we obtain
\begin{eqnarray}
S_k(w_1)S_l(w_2)=\left(\frac{w_1}{w_2}\right)^{A_{k,l}(0)}\frac{h_{k,l}\left(\frac{w_2}{w_1}\right)}{
h_{l,k}\left(\frac{w_1}{w_2}\right)}
S_l(w_2)S_k(w_1)~~~(1\leq k, l \leq L).
\nonumber
\end{eqnarray}
Using (\ref{thm:sj(m)}), (\ref{thm:Aij(0)}) and (\ref{thm:Aij(m)}), 
we obtain (\ref{eqn:screening}).
~~$\Box$\\

By direct calculation, we have the following lemma.
\begin{lem}
\label{lemma:3-12}
The determinants of $(A_{i,j}(m))_{i,j=1}^L$ in (\ref{thm:Aij(0)}) and (\ref{thm:Aij(m)}) are given by
\begin{eqnarray}
&&
\det\left(\left(A_{i,j}(m)\right)_{i,j=1}^L\right)=(-1)^L
\frac{\displaystyle
[D(0,L;\Phi)m]_x[(r-1)m]_x^{|\widehat{I}(\frac{1-r}{r})|-1}
[m]_x^{|\widehat{I}(-\frac{1}{r})|-1}
}{\displaystyle
[rm]_x^L \prod_{j=1}^L s_j(m)s_j(-m)}~~~(m\neq 0),
\nonumber
\\
&&
\det\left(\left(A_{i,j}(0)\right)_{i,j=1}^L\right)=
r^{-L} (r-1)^{|\widehat{I}(\frac{1-r}{r})|-1} D(0,L;\Phi).\nonumber
\end{eqnarray}
\end{lem}
Hence, the condition
$\det\left(\left(A_{i,j}(m)\right)_{i,j=1}^L\right)\neq 0$ 
$(m \in \mathbf{Z})$ is satisfied.

\section{Quadratic relation}
\label{Section4}
In this section, we introduce the higher $W$-currents $T_i(z)$
and obtain a set of quadratic relations of $T_i(z)$ 
for the deformed
$W$-superalgebra ${\cal W}_{q,t}\bigl(A(M,N)\bigr)$.
We show that these relations are
independent of the choice of Dynkin-diagrams.

\subsection{Quadratic relation}

We define the functions $\Delta_i(z)$ $(i=0,1,2, \ldots)$ as
\begin{eqnarray}
\Delta_i(z)=\frac{(1-x^{2r-i}z)(1-x^{-2r+i}z)}{(1-x^i z)(1-x^{-i}z)}.
\nonumber
\end{eqnarray}
We have
\begin{eqnarray}
&&
\Delta_i(z)-\Delta_i({z}^{-1})=\frac{[r]_x[r-i]_x}{[i]_x}(x-x^{-1})(\delta(x^{-i}z)-\delta(x^iz))~~(i=1,2,3,\ldots).
\nonumber
\end{eqnarray}

We define the structure functions $f_{i,j}(z;a)$ 
$(i, j=0,1,2,\ldots)$ as
\begin{eqnarray}
f_{i,j}(z;a)=\exp\left(-\sum_{m=1}^\infty \frac{1}{m}
\frac{[(r-1)m]_x[rm]_x
[{\rm Min}(i,j)m]_x[(a-{\rm Max}(i,j))m]_x
}{[m]_x[a m]_x}(x-x^{-1})^2
z^m
\right).
\label{def:fij}
\end{eqnarray}
In the case of $a=D(0,L; \Phi)$, 
the ratio of the structure function
\begin{eqnarray}
\frac{f_{1,1}(z^{-1};a)}{f_{1,1}(z;a)}=
\frac{\Theta_{x^{2a}}(x^2z, x^{-2r}z, x^{2r-2}z)}{
\Theta_{x^{2a}}(x^{-2}z,x^{2r}z, x^{-2r+2}z)}
\nonumber
\end{eqnarray}
coincides with those of (\ref{eqn:vertex}).

We introduce the higher $W$-currents $T_i(z)$ and give the quadratic relations.
From now on, we set $g=1$ in (\ref{thm:g}), but this is not an essential limitation.
Hereafter, we use the abbreviations
\begin{eqnarray}
c(r, x)=[r]_x[r-1]_x(x-x^{-1}),~~~d_j(r, x)=
\left\{\begin{array}{cc}
\displaystyle{\prod_{l=1}^j \frac{[r-l]_x}{[l]_x}}
&(j \geq 1),\\
1&(j=0).
\end{array}\right.
\nonumber
\end{eqnarray}

We introduce the $W$-currents $T_i(z)$ $(i=0,1,2, \ldots)$ as
\begin{align}
T_0(z)&=1,~~~
T_1(z)=\sum_{k \in \widehat{I}(\frac{1-r}{r})}
\Lambda_k(z)+d_1(r,x)
\sum_{k \in \widehat{I}(-\frac{1}{r})}\Lambda_k(z),
\notag
\\
T_i(z)&=\sum_{
(m_1,m_2,\ldots, m_{L+1}) \in \hat{N}(\Phi)
\atop{
m_1+m_2+\cdots+m_{L+1}=i}}
\prod_{k \in \hat{I}(-\frac{1}{r})}
d_{m_k}(r,x)~\Lambda_{m_1, m_2, \ldots, m_{L+1}}^{(i)}(z)~~(i=2,3,4,\ldots).
\label{def:Ti(z)}
\end{align}
Here we set
\begin{align}
\Lambda_{m_1,m_2,\cdots, m_{L+1}}^{(i)}(z)
&=:\prod_{k \in \widehat{I}(\frac{1-r}{r})}
\Lambda_k^{(m_k)}(x^{-i+1+2(m_1+\cdots+m_{k-1})}z)
\notag
\\
&\times 
\prod_{k \in \widehat{I}(-\frac{1}{r})}
\Lambda_{k}^{(m_k)}(x^{-i+1+2(m_1+\cdots+m_{k-1})}z):~~
{\rm for}~~(m_1,m_2, \ldots, m_{L+1}) \in \hat{N}(\Phi),
\notag
\end{align}
where
\begin{align}
&\Lambda_{k}^{(0)}(z)=1,~~~
\Lambda_{k}^{(m)}(z)=:\Lambda_{k}(z)
\Lambda_{k}(x^{2}z)\cdots \Lambda_{k}(x^{2m-2}z):,
\notag
\\
&
\hat{N}(\Phi)=
\left\{(m_1, m_2, \ldots, m_{L+1}) \in {\mathbf N}^{L+1}
\left|
0\leq m_k \leq 1~~{\rm if}~~k \in \hat{I}\left(\frac{1-r}{r}\right),~
m_k \geq 0~~{\rm if}~~k \in \hat{I}\left(-\frac{1}{r}\right)\right.
\right\}.
\label{def:N(Phi)}
\end{align}

We have $T_i(z)\neq 1$ $(i=1,2,3,\ldots)$ and $T_i(z) \neq T_j(z)$ $(i \neq j)$.

The following is {\bf the main theorem} of this paper.

\begin{thm}
\label{thm:4-1}~
For the deformed $W$-superalgebra 
${\cal W}_{q,t}\bigl(A(M,N)\bigr)$, 
the $W$-currents $T_i(z)$
satisfy the set of quadratic relations
\begin{align}
&
f_{i,j}\left(\frac{z_2}{z_1}; a\right)T_i(z_1)T_j(z_2)-
f_{j,i}\left(\frac{z_1}{z_2}; a\right)T_j(z_2)T_i(z_1)\notag\\
=&~c(r,x)\sum_{k=1}^i \prod_{l=1}^{k-1}\Delta_1(x^{2l+1})
\left(
\delta\left(\frac{x^{-j+i-2k}z_2}{z_1}\right)f_{i-k, j+k}(x^{j-i};a)T_{i-k}(x^{k}z_1)T_{j+k}(x^{-k}z_2)
\right.\notag
\\
-&\left.\delta\left(\frac{x^{j-i+2k}z_2}{z_1}\right)
f_{i-k, j+k}(x^{-j+i};a)T_{i-k}(x^{-k}z_1)T_{j+k}(x^kz_2)
\right)~~~(j \geq i \geq 1).
\label{thm:quadratic}
\end{align}
Here we use $f_{i,j}(z;a)$ in (\ref{def:fij}) with the specialization $a=D(0, L; \Phi)$.
The quadratic relations (\ref{thm:quadratic})
are independent of the choice of Dynkin-diagrams for the Lie superalgebra $A(M,N)$.
\end{thm}

In view of Theorem \ref{thm:4-1}, we arrive at the following definition. 
\begin{dfn}
\label{def:4-2}
Set $T_i(z)=\sum_{m \in {\mathbf Z}}T_i[m]z^{-m}$ $(i=1,2,3,\ldots)$.
The deformed $W$-superalgebra 
${\cal W}_{q, t}\bigl(A(M,N)\bigr)$
is an associative algebra over ${\mathbf C}$ with the generators $T_i[m]$ $(m \in {\mathbf Z}, i=1,2,3,\ldots)$
and the defining relations (\ref{thm:quadratic}).
\end{dfn}

\subsection{Proof of Theorem \ref{thm:4-1}}

\begin{prop}
\label{prop:4-3}
The $\Lambda_i(z)$'s satisfy
\begin{eqnarray}
&&
~f_{1,1}\left(\frac{z_2}{z_1};a\right)
\Lambda_k(z_1)\Lambda_l(z_2)=
\Delta_1\left(\frac{x^{-1}z_2}{z_1}\right)
:\Lambda_k(z_1)\Lambda_l(z_2):,
\nonumber
\\
&&
f_{1,1}\left(\frac{z_2}{z_1};a\right)
\Lambda_l(z_1)\Lambda_k(z_2)=
\Delta_1\left(\frac{x z_2}{z_1}\right)
:\Lambda_l(z_1)\Lambda_k(z_2):
~~~(1\leq k<l\leq L+1),
\nonumber
\\
&&
f_{1,1}\left(\frac{z_2}{z_1};a\right)
\Lambda_i(z_1)\Lambda_i(z_2)=
:\Lambda_i(z_1)\Lambda_i(z_2):~~~
{\rm if}~~i \in \hat{I}\left(\frac{1-r}{r}\right),
\nonumber
\\
&&
f_{1,1}\left(\frac{z_2}{z_1};a\right)
\Lambda_i(z_1)\Lambda_i(z_2)=
\Delta_2\left(\frac{z_2}{z_1}\right):
\Lambda_i(z_1)\Lambda_i(z_2):~~~
{\rm if}~~i \in \hat{I}\left(-\frac{1}{r}\right)
\label{normal:lambda(z)}
\end{eqnarray}
where we set $a=D(0,L ; \Phi)$.
\end{prop}
{\it Proof of Proposition \ref{prop:4-3}.}~
Substituting 
(\ref{thm:pij}) and (\ref{thm:lambda(m)}) into 
$\varphi_{\Lambda_k, \Lambda_l}(z_1,z_2)$ in (\ref{normal4}), we obtain (\ref{normal:lambda(z)}).~~
$\Box$

~\\
{\it Proof of Proposition \ref{prop:3-2}.}
\label{proof:prop:3-2}~
Using $\varphi_{\Lambda_k, \Lambda_l}(z_1,z_2)$ in (\ref{normal4}), we obtain
\begin{eqnarray}
\Lambda_k(z_1)\Lambda_l(z_2)=\frac{\varphi_{\Lambda_k, \Lambda_l}\left(z_1, z_2\right)}{
\varphi_{\Lambda_l, \Lambda_k}\left(z_2, z_1\right)}
\Lambda_l(z_2)\Lambda_k(z_1)~~~(1\leq k, l \leq L+1).
\nonumber
\end{eqnarray}
Using the explicit formulae of $\varphi_{\Lambda_k, \Lambda_l}(z_1,z_2)$,
we obtain (\ref{eqn:vertex}).
~~$\Box$\\

\begin{lem}
\label{lem:4-4}
The $D(0,L; \Phi)$ given in (\ref{def:D(k,l)}) 
is independent of the choice of the 
Dynkin-diagrams for the Lie superalgebra $A(M,N)$.
\begin{eqnarray}
D(0,L;\Phi)=D(0,L; r_{\alpha_i}(\Phi))~~~
(\alpha_i \in \Pi).\label{eqn:D-invariance}
\end{eqnarray}
Here $\Pi$ is a fundamental system.
\end{lem}
{\it Proof of Lemma \ref{lem:4-4}.}~
We show (\ref{eqn:D-invariance}) by checking all cases.
We set the Dynkin-diagrams $\Phi_j$ $(1\leq j \leq 8)$ as follows.
Let $K=K(\Phi_j)$ the number of odd isotropic roots $(\alpha_i, \alpha_i)=0$
in the Dynkin-diagram $\Phi_j$.
We set
\begin{center}
\input{Dynkin16.tex}
~~~~~~~~~~~~~~~
\input{Dynkin17.tex}
\end{center}
\begin{center}
~~~~~~~~~~~
\input{Dynkin18.tex}
~~~
\input{Dynkin19.tex}
\end{center}
For $2\leq i \leq L-1$, we set
\begin{center}
\input{Dynkin20.tex}
\end{center}

\begin{center}
\input{Dynkin21.tex}
\end{center}

\begin{center}
\input{Dynkin22.tex}
\end{center}

\begin{center}
\input{Dynkin23.tex}
\end{center}
We have
$r_{\alpha_1}(\Phi_1)=\Phi_2$,
$r_{\alpha_L}(\Phi_3)=\Phi_4$,
$r_{\alpha_i}(\Phi_5)=\Phi_6$,
and $r_{\alpha_i}(\Phi_7)=\Phi_8$.

The affinized Dynkin-diagrams
$\widehat{\Phi}_j$ from $\Phi_j$ $(1\leq j \leq 8)$
are given as 
$
\widehat{\Phi}_j=\left\{\begin{array}{cc}
\widehat{\Phi}_{j,1}&{\rm if}~K(\Phi_j)={\rm even},\\
\widehat{\Phi}_{j,2}&{\rm if}~K(\Phi_j)={\rm odd}.
\end{array}\right.$

\vspace{6mm}
\begin{center}
\input{Dynkin24.tex}
~~~~~
\input{Dynkin25.tex}
~~~~~~~~~\\
\vspace{8mm}
\input{Dynkin26.tex}
~~~~~
\input{Dynkin27.tex}
~~~~~~~~~\\
\vspace{8mm}
\input{Dynkin28.tex}
\input{Dynkin29.tex}
~~~~~~~~~~~~\\
\vspace{6mm}
\input{Dynkin30.tex}
\input{Dynkin31.tex}
~~~~~~~~~~~~
\end{center}
\vspace{1mm}
\begin{center}
\input{Dynkin32.tex}
\\
\vspace{6mm}
\input{Dynkin33.tex}
\end{center}
\vspace{7mm}
\begin{center}
\input{Dynkin34.tex}
\\
\vspace{6mm}
\input{Dynkin35.tex}
\end{center}
\vspace{7mm}
\begin{center}
\input{Dynkin36.tex}
\\
\vspace{6mm}
\input{Dynkin37.tex}
\end{center}
\vspace{7mm}
\begin{center}
\input{Dynkin38.tex}
\\
\vspace{6mm}
\input{Dynkin39.tex}
\end{center}
\vspace{3mm}
The values of $A_{j-1,j}(0)$ are
written beside the line segment connecting $\alpha_{j-1}$
and $\alpha_j$.
We have
$$|\widehat{I}(1,L+1;\delta,\widehat{\Phi}_{2j-1,1})|=
|\widehat{I}(1,L+1;\delta,\widehat{\Phi}_{2j,2})|,
~~~|\widehat{I}(1,L+1;\delta,\widehat{\Phi}_{2j-1,2})|=
|\widehat{I}(1,L+1;\delta,\widehat{\Phi}_{2j,1})|~~(1\leq j \leq 2),
$$ 
$$
|\widehat{I}(1,L+1;\delta,\widehat{\Phi}_{2j-1,1})|=
|\widehat{I}(1,L+1;\delta,\widehat{\Phi}_{2j,1})|,
~~~|\widehat{I}(1,L+1;\delta,\widehat{\Phi}_{2j-1,2})|=
|\widehat{I}(1,L+1;\delta,\widehat{\Phi}_{2j,2})|~~(3\leq j \leq 4),
$$ 
where $\delta=\frac{1-r}{r}$ or $-\frac{1}{r}$.
Hence we have
$$
D(0,L;\Phi_{2j-1})=D(0,L;\Phi_{2j})~~~(1\leq j \leq 4).$$
In other words, we have
\begin{eqnarray}
&&
D(0,L;\Phi_j)=D(0,L;r_{\alpha_1}(\Phi_j))~~(1\leq j \leq 2),\nonumber
\\
&&
D(0,L;\Phi_j)=D(0,L;r_{\alpha_L}(\Phi_j))~~(3\leq j \leq 4),
\nonumber\\
&&
D(0,L;\Phi_j)=D(0,L;r_{\alpha_i}(\Phi_j))~~(5 \leq j \leq 8,
2\leq i \leq L-1).\nonumber
\end{eqnarray}
Now we have proved (\ref{eqn:D-invariance}).~~$\Box$

\begin{lem}
\label{lem:4-5}
$\Delta_i(z)$ and $f_{i,j}(z;a)$ satisfy the following fusion relations.
\begin{eqnarray}
&&
f_{i,j}(z;a)=f_{j,i}(z;a)=
\prod_{k=1}^i f_{1,j}(z^{-i-1+2k}z; a)~~~(1\leq i \leq j),
\label{eqn:fusion0}
\\
&&
f_{1,i}(z;a)=
\left(
\prod_{k=1}^{i-1}\Delta_1(x^{-i+2k}z)
\right)^{-1}
\prod_{k=1}^i f_{1,1}(x^{-i-1+2k}z;a)~~~
(i \geq 2),
\label{eqn:fusion5}
\\
&&
f_{1,i}(z;a)f_{j,i}(x^{\pm (j+1)}z;a)=\left\{\begin{array}{cc}
f_{j+1,i}(x^{\pm j}z;a)\Delta_1(x^{\pm i}w)&(1\leq i \leq j),\\
f_{j+1,i}(x^{\pm j}z;a)&(1\leq j<i),
\end{array}\right.
\label{eqn:fusion1}\\
&&
f_{1,i}(z;a)f_{1,j}(x^{\pm (i+j)}z;a)=
f_{1,i+j}(x^{\pm j}z;a)\Delta_1(x^{\pm i}z)~~(i,j \geq 1),
\label{eqn:fusion2}\\
&&
f_{1,i}(z;a)
f_{1,j}(x^{\pm (i-j-2k)}z;a)=
f_{1,i-k}(x^{\mp k}z;a)f_{1,j+k}(x^{\pm (i-j-k)}z;a)~~(i,j,i-k,j+k \geq 1),
\label{eqn:fusion3}
\\
&&
\Delta_{i+1}(z)=
\left(\prod_{k=1}^{i-1} \Delta_{1}(x^{-i+2k}z)
\right)^{-1}
\prod_{k=1}^{i}\Delta_2(x^{-i-1+2k}z)~~~
(i \geq 2).
\label{eqn:fusion4}
\end{eqnarray}
\end{lem}
{\it Proof of Lemma \ref{lem:4-5}.}~
We obtain (\ref{eqn:fusion0}) and (\ref{eqn:fusion4})
by straightforward calculation from the definitions.
We show (\ref{eqn:fusion5}) here.
From definitions, we have
\begin{align}
&
\left(\prod_{k=1}^{i-1}\Delta_1(x^{-i+2k}z)\right)^{-1}
\prod_{k=1}^i f_{1,1}(x^{-i-1+2k}z;a)\nonumber
\\
=&
\exp\left(
-\sum_{m=1}^\infty
\frac{1}{m}
\frac{[rm]_x[(r-1)m]_x}{[a m]_x}
(x-x^{-1})^2
\left([(a-1)m]_x\sum_{k=1}^i x^{(-i+2k-1)m}-
[a m]_x\sum_{k=1}^{i-1}
x^{(-i+2k)m}
\right)z^m
\right).\nonumber
\end{align}
Using the relation
\begin{eqnarray}
[(a-1)m]_x\sum_{k=1}^i x^{(-i+2k-1)m}-
[a m]_x\sum_{k=1}^{i-1}
x^{(-i+2k)m}=[(a-i)m]_x,
\nonumber
\end{eqnarray}
we have $f_{1,i}(z;a)$.
Using (\ref{eqn:fusion0}) and (\ref{eqn:fusion5}),
we obtain the relations (\ref{eqn:fusion1}), (\ref{eqn:fusion2}), and (\ref{eqn:fusion3}).~~~
$\Box$

The following relations
(\ref{eqn:fusion6}), (\ref{eqn:TiTj-commute}), and (\ref{def:quadratic'}) give special cases of (\ref{thm:quadratic}).
\begin{lem}
\label{lem:4-6}~
The $T_i(z)$'s satisfy the fusion relation
\begin{eqnarray}
&&
\lim_{z_1 \to x^{\pm (i+j)}z_2}
\left(1-x^{\pm(i+j)}\frac{z_2}{z_1}\right)
f_{i,j}\left(\frac{z_2}{z_1};a\right)T_i(z_1)T_j(z_2)\nonumber
\\
&&=\mp c(r,x) \prod_{k=1}^{{\rm Min}(i,j)-1}\Delta_1(x^{2k+1})T_{i+j}(x^{\pm i}z_2)~~~(i,j \geq 1).
\label{eqn:fusion6}
\end{eqnarray}
Here we set $a=D(0,L;\Phi)$.
\end{lem}
{\it Proof of Lemma \ref{lem:4-6}.}~
Summing up the relations 
(\ref{appendixA:1})--(\ref{appendixA:4}) in Appendix \ref{Appendix:fusion}
gives (\ref{eqn:fusion6}).~~$\Box$

\begin{lem}
\label{lem:4-7}~
The $T_i(z)$'s satisfy the exchange relation as 
meromorphic functions
\begin{eqnarray}
f_{i,j}\left(\frac{z_2}{z_1};a\right)T_i(z_1)T_j(z_2)
=f_{j,i}\left(\frac{z_1}{z_2};a\right)T_j(z_2)T_i(z_1)~~~(j \geq i \geq 1).
\label{eqn:TiTj-commute}
\end{eqnarray}
Both sides are regular except for poles at
$z_2/z_1 \neq x^{\pm(j-i+2k)}$ $(1 \leq k \leq i)$.
Here we set $a=D(0,L;\Phi)$.
\end{lem}
{\it Proof of Lemma \ref{lem:4-7}.}~Using the commutation relation (\ref{eqn:vertex}) repeatedly, 
(\ref{eqn:TiTj-commute}) is obtained except for poles in both sides.
Using Proposition \ref{prop:4-3},
we identify the pole position
as $z_2/z_1=x^{\pm(j-i+2k)}$ $(1 \leq k \leq i)$.
~~$\Box$

\begin{lem}
\label{lem:4-8}~
The $T_i(z)$'s satisfy the quadratic relations
\begin{align}
&
f_{1,i}\left(\frac{z_2}{z_1};a\right)T_1(z_1)T_i(z_2)
-f_{i,1}\left(\frac{z_1}{z_2};a\right)T_i(z_2)T_1(z_1)\nonumber\\
=&~c(r,x)\left(\delta\left(\frac{x^{-i-1}z_2}{z_1}\right)T_{i+1}(x^{-1}z_2)-\delta\left(\frac{x^{i+1}z_2}{z_1}\right)T_{i+1}(xz_2)\right)
~~(i \geq 1).
\label{def:quadratic'}
\end{align}
Here we set $a=D(0,L;\Phi)$.
\end{lem}
{\it Proof of Lemma \ref{lem:4-8}}.~
Summing up the relations 
(\ref{appendixB:2})--(\ref{appendixB:6}) in Appendix \ref{Appendix:exchange} 
gives (\ref{def:quadratic'}).~~
$\Box$

~\\
~{\it Proof of Theorem \ref{thm:4-1}}.~
We prove Theorem \ref{thm:4-1} by induction.
Lemma \ref{lem:4-8} is the basis of induction for the proof.
In what follows we set $a=D(0,L;\Phi)$.

We define ${\rm LHS}_{i,j}$ and ${\rm RHS}_{i,j}(k)$ with $(1\leq k \leq i \leq j)$ as
\begin{align}
{\rm LHS}_{i,j}&=
f_{i,j}\left(\frac{z_2}{z_1};a\right)T_i(z_1)T_j(z_2)
-f_{j,i}\left(\frac{z_1}{z_2};a\right)T_j(z_2)T_i(z_1),
\notag
\\
{\rm RHS}_{i,j}(k)&=
c(r,x)
\prod_{l=1}^{k-1}\Delta_1(x^{2l+1})
\left(\delta\left(\frac{x^{-j+i-2k}z_2}{z_1}\right)
f_{i-k,j+k}(x^{j-i};a)T_{i-k}(x^kz_1)T_{j+k}(x^{-k}z_2)\right.
\notag\\
&-\left.\delta
\left(\frac{x^{j-i+2k}z_2}{z_1}\right)f_{i-k,j+k}(x^{-j+i};a)T_{i-k}(x^{-k}z_1)T_{j+k}(x^{k}z_2)
\right)~~
(1\leq k \leq i-1),
\notag
\\
{\rm RHS}_{i,j}(i)&=
c(r,x)
\prod_{l=1}^{i-1}\Delta_1(x^{2l+1})
\left(\delta\left(\frac{x^{-j-i}z_2}{z_1}\right)T_{j+i}(x^{-i}z_2)-
\delta\left(\frac{x^{j+i}z_2}{z_1}\right)
T_{j+i}(x^{i}z_2)\right).
\notag
\end{align}
We prove the following relation by induction on 
$i$ $(1\leq i \leq j)$.
\begin{eqnarray}
{\rm LHS}_{i,j}=\sum_{k=1}^i {\rm RHS}_{i,j}(k).
\label{induction4}
\end{eqnarray}
The starting point of $i=1 \leq j$ was previously proven in 
Lemma \ref{lem:4-8}.
We assume that the relation (\ref{induction4}) holds for $i$ $(1\leq i<j)$, and we show ${\rm LHS}_{i+1,j}=\sum_{k=1}^{i+1} {\rm RHS}_{i+1,j}(k)$ from this assumption.
Multiplying ${\rm LHS}_{i,j}$ by 
$f_{1,i}\left(z_1/z_3;a\right)f_{1,j}\left(z_2/z_3;a\right)T_1(z_3)$ on the left and using the quadratic relation 
$f_{1,j}\left(z_2/z_3;a\right)T_1(z_3)T_j(z_2)=
f_{j,1}\left(z_3/z_2;a\right)T_j(z_2)T_1(z_3)+\cdots$, along with the fusion relation (\ref{eqn:fusion1}) gives 
\begin{align}
&f_{1,j}\left(\frac{z_2}{z_3};a\right)
f_{i,j}\left(\frac{z_2}{z_1};a\right)
f_{1,i}\left(\frac{z_1}{z_3};a\right)T_1(z_3)T_i(z_1)T_j(z_2)\notag\\
-&f_{j,1}\left(\frac{z_3}{z_2};a\right)
f_{j,i}\left(\frac{z_1}{z_2};a\right)T_j(z_2)
f_{1,i}\left(\frac{z_1}{z_3};a\right)T_1(z_3)T_i(z_1)\notag\\
-&~c(r,x)\delta\left(\frac{x^{-j-1}z_2}{z_3}\right)
\Delta_1\left(\frac{x^{-i}z_1}{z_3}\right)f_{j+1,i}\left(\frac{x^{-j}z_1}{z_3};a\right)T_{j+1}(x^jz_3)T_i(z_1)\notag
\\
+&~c(r,x)\delta\left(\frac{x^{j+1}z_2}{z_3}\right)
\Delta_1\left(\frac{x^{i}z_1}{z_3}\right)
f_{j+1,i}\left(\frac{x^{j}z_1}{z_3};a\right)T_{j+1}(x^{-j}z_3)T_i(z_1).
\label{proof1}
\end{align}
Taking the limit $z_3\to x^{-i-1}z_1$ of (\ref{proof1}) multiplied by $c(r,x)^{-1} \left(1-x^{-i-1}z_1/z_3\right)$ and using the fusion relation (\ref{eqn:fusion6}) along with the relation $\lim_{z_3 \to x^{-i-1}z_1}\left(1-x^{-i-1}z_1/z_3\right)
\Delta_1\left(x^{-i}z_1/z_3\right)=c(r,x)$ gives
\begin{align}
&f_{1,j}\left(\frac{x^{i+1}z_2}{z_1};a\right)
f_{i,j}\left(\frac{z_2}{z_1};a\right)
T_{i+1}(x^{-1}z_1)T_j(z_2)-
f_{j,1}\left(\frac{x^{-i-1}z_1}{z_2};a\right)
f_{j,i}\left(\frac{z_1}{z_2};a\right)
T_j(z_2)T_{i+1}(x^{-1}z_1)\notag\\
-&~c(r,x)\delta\left(\frac{x^{i-j}z_2}{z_1}\right)
f_{j+1,i}(x^{i-j+1};a)T_{j+1}(x^{j-i-1}z_1)T_i(z_1)\notag\\
-&~c(r,x)\delta\left(\frac{x^{i+j+2}z_2}{z_1}\right)
\prod_{l=1}^i \Delta_1(x^{2l+1})T_{i+j+1}(x^{i+1}z_2).
\notag
\end{align}
Using the fusion relation (\ref{eqn:fusion1}) and 
$f_{j+1,i}(x^{i-j+1};a)T_{j+1}(x^{j-i-1}z_1)T_i(z_1)=
f_{i,j+1}(x^{j-i-1};a)T_i(z_1)T_{j+1}(x^{j-i-1}z_1)$ in (\ref{eqn:TiTj-commute})
gives
\begin{align}
&
f_{i+1,j}\left(\frac{xz_2}{z_1};a\right)
T_{i+1}(x^{-1}z_1)T_j(z_2)-
f_{j,i+1}\left(\frac{x^{-1}z_1}{z_2};a\right)
T_j(z_2)T_{i+1}(x^{-1}z_1)
\notag
\\
-&~c(r,x)\delta\left(\frac{x^{i-j}z_2}{z_1}\right)
f_{i,j+1}(x^{-i+j+1};a)T_i(z_1)T_{j+1}(x^{-1}z_2)
\notag\\
+&~c(r,x)\delta\left(\frac{x^{i+j+2}z_2}{z_1}\right)\prod_{l=1}^{i}\Delta_1(x^{2l+1})T_{i+j+1}(x^{i+1}z_2).
\label{induction1}
\end{align}
Multiplying ${\rm RHS}_{i,j}(i)$ by $f_{1,i}\left(z_1/z_3;a\right)
f_{1,j}\left(z_2/z_3;a\right)T_1(z_3)$ from the left and using the fusion relation (\ref{eqn:fusion2}) gives
\begin{align}
&c(r,x)\prod_{l=1}^{i-1}\Delta_1(x^{2l+1})
\left(\delta\left(\frac{x^{-i-j}z_2}{z_1}\right)
f_{1,i+1}\left(\frac{x^jz_1}{z_3};a\right)\Delta_1\left(\frac{x^i z_1}{z_3}\right)T_1(z_3)T_{i+j}(x^jz_1)\right.
\notag\\
-&~\left.
\delta\left(\frac{x^{i+j}z_2}{z_1}\right)
f_{1,i+1}\left(\frac{x^{-j}z_1}{z_3};a\right)
\Delta_1\left(\frac{x^{-i} z_1}{z_3}\right)T_1(z_3)T_{i+j}(x^{-j}z_1)
\right).
\label{proof2}
\end{align}
Taking the limit $z_3\to x^{-i-1}z_1$ of (\ref{proof2}) multiplied by $c(r,x)^{-1}\left(1-x^{-i-1}z_1/z_3\right)$ 
and using the fusion relation (\ref{eqn:fusion6}) along with the relation $\lim_{z_3 \to x^{-i-1}z_1}\left(1-x^{-i-1}z_1/z_3\right)
\Delta_1\left(x^{-i}z_1/z_3\right)=c(r,x)$ gives
\begin{align}
&
c(r,x)\delta\left(\frac{x^{-i-j}z_2}{z_1}\right)\prod_{l=1}^i \Delta_1(x^{2l+1})T_{i+j+1}(x^{-i-1}z_2)\notag\\
-~&
c(r,x)\delta\left(\frac{x^{i+j}z_2}{z_1}\right)
\prod_{l=1}^{i-1} \Delta_1(x^{2l+1})
f_{1,i+j}(x^{i-j+1};a)T_1(x^{-i-1}z_1)
T_{i+j}(x^{i}z_2).\label{induction2}
\end{align}
Multiplying ${\rm RHS}_{i,j}(k)$ $(1\leq k \leq i-1)$ by $f_{1,i}\left(z_1/z_3;a\right)f_{1,j}\left(z_2/z_3;a\right)T_1(z_3)$ 
from the left and using the fusion relation (\ref{eqn:fusion3}) along with $f_{i-k,j+k}(x^{j-i};a)T_{i-k}(x^kz_1)T_{j+k}(x^{j-i+k}z_1)=f_{j+k,i-k}(x^{i-j};a)T_{j+k}(x^{j-i+k}z_1)T_{i-k}(x^kz_1)$ in (\ref{eqn:TiTj-commute}) 
gives
\begin{align}
&
c(r,x)\prod_{l=1}^{k-1}\Delta_1(x^{2l+1})
\label{proof3}
\\
\times&\left(
\delta\left(\frac{x^{-j+i-2k}z_2}{z_1}\right)
f_{1,i-k}\left(\frac{x^kz_1}{z_3};a\right)f_{j+k,i-k}(x^{i-j};a)f_{1,j+k}\left(\frac{x^{-i+j+k}z_1}{z_3};a\right)T_1(z_3)T_{j+k}(x^{j-i+k}z_1)T_{i-k}(x^kz_1)
\right.\notag\\
&-\left.
\delta\left(\frac{x^{j-i+2k}z_2}{z_1}\right)
f_{1,i-k}\left(\frac{x^{-k}z_1}{z_3};a\right)f_{i-k,j+k}(x^{i-j};a)
f_{1,j+k}\left(\frac{x^{i-j-k}z_1}{z_3};a\right)T_1(z_3)T_{i-k}(x^{-k}z_1)T_{j+k}(x^kz_2)
\right).\nonumber
\end{align}
Taking the limit $z_3\to x^{-i-1}z_1$ of (\ref{proof3}) multiplied by $c(r,x)^{-1}\left(1-x^{-i-1}z_1/z_3\right)$ 
and using the fusion relations (\ref{eqn:fusion1}) and (\ref{eqn:fusion6}) along with
$$
f_{i-k+1,j+k}(x^{i-j+1};a)T_{i-k+1}(x^{-k-1}z_1)T_{j+k}(x^{-j+i-k}z_1)
=f_{j+k,i-k+1}(x^{j-i-1};a)T_{j+k}(x^{-j+i-k}z_1)T_{i-k+1}(x^{-k-1}z_1)
$$ in (\ref{eqn:TiTj-commute})
gives
\begin{align}
&
c(r,x)\prod_{l=1}^k \Delta_1(x^{2l+1})
\delta\left(\frac{x^{-j+i-2k}z_2}{z_1}\right)
f_{j+k-1,i-k}(x^{i-j+1};a)
T_{i-k}(x^kz_1)T_{j+k+1}(x^{-k-1}z_2)\nonumber
\\
-~&
c(r,x)
\prod_{l=1}^{k-1}\Delta_1(x^{2l+1})
\delta\left(\frac{x^{j-i+2k}z_2}{z_1}\right)
f_{i-k+1,j+k}(x^{i-j+1};a)T_{i-k+1}(x^{-k-1}z_1)T_{j+k}(x^kz_2).
\label{induction3}
\end{align}
Summing (\ref{induction1}), (\ref{induction2}), and (\ref{induction3}) for $1\leq k \leq i-1$ and shifting the variable $z_1 \mapsto xz_1$ gives ${\rm LHS}_{i+1,j}=\sum_{k=1}^{i+1}{\rm RHS}_{i+1,j}(k)$.
By induction on $i$, we have shown the quadratic relation (\ref{thm:quadratic}).

The quadratic relations (\ref{thm:quadratic}) are independent of 
the choice of Dynkin-diagrams for the Lie superalgebra $A(M,N)$,
because $a=D(0,L;\Phi)$ is independent of 
the choice of Dynkin-diagrams. See Lemma \ref{lem:4-4}.
~~$\Box$

\subsection{Classical limit}

The deformed $W$-algebra 
${\cal W}_{q, t}\bigl({\mathfrak g}\bigr)$ includes the $q$-Poisson $W$-algebra as a special case.
As an application of the quadratic relations (\ref{thm:quadratic}),
we obtain the $q$-Poisson $W$-algebra \cite{Frenkel-Reshetikhin, Frenkel-Reshetikhin1,Frenkel-Reshetikhin-Semenov}. 
We study ${\cal W}_{q,t}(A(M,N))$ 
$(M \geq N \geq 0, M+N \geq 1)$.
We set parameters $q=x^{2r}$ and $\beta=(r-1)/r$.
We define the $q$-Poisson bracket $\{,\}$ by taking the classical limit $\beta \to 0$ with $q$ fixed as
\begin{eqnarray}
\{T_i^{{\rm PB}}[m], T_j^{{\rm PB}}[n]\}=-\lim_{\beta \to 0}
\frac{1}{\beta \log q}[T_i[m],T_j[n]].\nonumber
\end{eqnarray}
Here, we set $T_i^{PB}[m]$ as
$T_i(z)=\sum_{m \in {\mathbf Z}}T_i[m]z^{-m} \longrightarrow T_i^{PB}(z)=\sum_{m \in {\mathbf Z}}T_i^{PB}[m]z^{-m}$
$(\beta \to 0,~q~\rm{fixed})$.
The $\beta$-expansions of the structure functions are given as 
\begin{eqnarray}
&&f_{i,j}(z;a)=1+\beta \log q \sum_{m=1}^\infty
\frac{\left[\frac{1}{2} {\rm Min}(i,j)m \right]_q 
\left[\left(\frac{1}{2}({\rm Max}(i,j)-M-1)\right)m\right]_q}{[\frac{1}{2}(M+1)m]_q} (q-q^{-1})+O(\beta^2)~~(i, j \geq 1),
\nonumber\\
&&c(r,x)=-\beta \log q+O(\beta^2),
\nonumber
\end{eqnarray}
where $a=D(0,M+N+1;\Phi)=(N+1)r+M-N$.
\begin{prop}
\label{prop:4-9}
For the $q$-Poisson $W$-superalgebra for $A(M,N)$
$(M \geq N \geq 0, M+N \geq 1)$, the generating functions
$T_i^{PB}(z)$ satisfy
\begin{align}
\{T_i^{PB}(z_1),T_j^{PB}(z_2)\}
&=(q-q^{-1})C_{i,j}\left(\frac{z_2}{z_1}\right)T_i^{PB}(z_1)T_j^{PB}(z_2)
\notag
\\
&+\sum_{k=1}^i \delta\left(\frac{q^{\frac{-j+i}{2}-k}z_2}{z_1}\right)T_{i-k}^{PB}(q^{\frac{k}{2}}z_1)T_{j+k}^{PB}(q^{-\frac{k}{2}}z_2)
\notag
\\
&-\sum_{k=1}^i \delta\left(\frac{q^{\frac{j-i}{2}+k}z_2}{z_1}\right)T_{i-k}^{PB}(q^{-\frac{k}{2}}z_1)T_{j+k}^{PB}(q^{\frac{k}{2}}z_2)~~(1\leq i \leq j).
\notag
\end{align}
Here we set the structure functions $C_{i,j}(z)$ $(i,j\geq 1)$ as
\begin{eqnarray}
C_{i,j}(z)=\sum_{m \in {\mathbf Z}}
\frac{\left[\frac{1}{2} {\rm Min}(i,j)m \right]_q \left[\frac{1}{2}({\rm Max}(i,j)-M-1)m\right]_q}{[\frac{1}{2}(M+1)m]_q} z^m~~(i,j \geq 1).
\nonumber
\end{eqnarray}
\end{prop}
The structure functions satisfy
$C_{i,M+1}(z)=C_{M+1,i}(z)=0$ $(1\leq i \leq M+1)$.

\section{Conclusion and Discussion}
\label{Section5}

In this paper,
we found the free field construction 
of the basic $W$-current $T_1(z)$ (See (\ref{thm:lambda(0)}) 
and (\ref{thm:lambda(m)}))
and the screening currents $S_j(w)$ (See (\ref{thm:sj(m)}))
for the deformed $W$-superalgebra 
${\cal W}_{q,t}\bigl(A(M,N)\bigr)$.
Using the free field construction,
we introduced the higher $W$-currents $T_i(z)$ (See (\ref{def:Ti(z)})) and obtained a closed set of quadratic relations among them (See (\ref{thm:quadratic})).
These relations are independent of the choice of Dynkin-diagrams
for the Lie superalgebra $A(M,N)$.

Recently, Feigin, Jimbo, Mukhin, and Vilkoviskiy
\cite{Feigin-Jimbo-Mukhin-Vilkoviskiy}
introduced the free field construction of
the basic $W$-current and the screening currents
in types $A, B, C, D$ including twisted and supersymmetric cases
in terms of the quantum toroidal algebras.
Their motivation is to understand
a commutative family of integrals of motion associated with affine Dynkin-diagrams \cite{Feigin-Kojima-Shiraishi-Watanabe,
Kojima-Shiraishi}.
In the case of type $A$, their basic $W$-current $T_1(z)$
satisfies
\begin{eqnarray}
T_1(z_1)T_1(z_2)=\frac{\Theta_{\mu}\left(
q_1\frac{z_2}{z_1},~q_2\frac{z_2}{z_1},~q_3\frac{z_2}{z_1}\right)}{
\Theta_{\mu}\left(
q_1^{-1}\frac{z_2}{z_1},~q_2^{-1}\frac{z_2}{z_1},~q_3^{-1}\frac{z_2}{z_1}
\right)}T_1(z_2)T_1(z_1)~~~(q_1q_2q_3=1)
\label{toroidal:vertex}
\end{eqnarray}
in the sense of analytic continuation.
Upon the specialization $q_1=x^2, q_2=x^{-2r}, q_3=x^{2r-2}, \mu=x^{2a}$ $(a=D(0,L;\Phi))$, their commutation relation
(\ref{toroidal:vertex}) coincides with those of this paper
(See (\ref{eqn:vertex})).
In the case of $\mathfrak{sl}(N)$,
their basic $W$-current $T_1(z)$
coincides with those of \cite{Feigin-Kojima-Shiraishi-Watanabe,
Kojima-Shiraishi}, which gives a one-parameter deformation of
${\cal W}_{q,t}(\mathfrak{sl}(N))$
in Ref.\cite{Awata-Kubo-Odake-Shiraishi, Feigin-Frenkel}.
In the case of $A(M,N)$,
their basic $W$-current $T_1(z)$
gives a one-parameter deformation of
those of ${\cal W}_{q,t}\bigl(A(M,N)\bigr)$
in this paper.

It is still an open problem to find quadratic relations of the deformed $W$-algebra ${\cal W}_{q,t}({\mathfrak g})$,
except for ${\mathfrak g}=\mathfrak{sl}(N)$, $A_2^{(2)}$, and 
$A(M,N)$.
It seems to be possible to extend Ding-Feigin's construction to other Lie superalgebras and obtain
their quadratic relations.

~\\
{\bf ACKNOWLEDGMENTS}

The author would like to thank Professor Michio Jimbo very much
for carefully reading the manuscript and for giving lots of useful advice.
This work is supported by the Grant-in-Aid for Scientific Research {\bf C} (26400105) from the Japan Society for the Promotion of Science.

\begin{appendix}
\renewcommand{\theequation}{\Alph{section} \arabic{equation}}
\setcounter{equation}{0}

\section{Fusion relation}
\label{Appendix:fusion}
In this appendix we summarize the fusion relations of
$\Lambda_i(z)$.
We use the abbreviation
\begin{eqnarray}
F_{i,j}^{(\pm)}(z;a)=(1-x^{\pm (i+j)}z)f_{i,j}(z;a)~~~(a=D(0,L;\Phi)).\nonumber
\end{eqnarray}
For $(m_1,m_2,\ldots,m_{L+1}),
(n_1, n_2,\ldots, n_{L+1}) \in \hat{N}(\Phi)$
defined in (\ref{def:N(Phi)}),
we set
$i=m_1+m_2+\cdots+m_{L+1}$ and
$j=n_1+n_2+\cdots+n_{L+1}$.\\
$\bullet$~If ${\rm Max}\{1\leq k \leq L+1|n_k\neq 0\}
<{\rm Min}\{1\leq k \leq L+1|m_k \neq 0\}$ holds, we have
\begin{align}
&
\lim_{z_1 \to x^{i+j}z_2}
F_{i,j}^{(+)}\left(\frac{z_2}{z_1};a\right)
\Lambda_{m_1, m_2, \ldots,m_{L+1}}^{(i)}(z_1)
\Lambda_{n_1, n_2, \ldots, n_{L+1}}^{(j)}(z_2)\notag
\\
&=
-c(r,x)\prod_{l=1}^{{\rm Min}(i,j)-1}\Delta_1(x^{2l+1})
\Lambda_{m_1+n_1, m_2+n_2, \ldots, m_{L+1}+n_{L+1}}^{(i+j)}(x^{i}z_2).
\label{appendixA:1}
\end{align}
$\bullet$~If ${\rm Max}\{1\leq k \leq L+1|m_k\neq 0\}
<{\rm Min}\{1\leq k \leq L+1|n_k \neq 0\}$ holds, we have
\begin{align}
&
\lim_{z_1 \to x^{-(i+j)}z_2}
F_{i,j}^{(-)}\left(\frac{z_2}{z_1};a\right)
\Lambda_{m_1, m_2, \ldots,m_{L+1}}^{(i)}(z_1)
\Lambda_{n_1, n_2, \ldots, n_{L+1}}^{(j)}(z_2)\notag
\\
&=
c(r,x)\prod_{l=1}^{{\rm Min}(i,j)-1}\Delta_1(x^{2l+1})
\Lambda_{m_1+n_1, m_2+n_2, \ldots, m_{L+1}+n_{L+1}}^{(i+j)}(x^{-i}z_2).
\label{appendixA:2}
\end{align}
$\bullet$~If $l$ satisfies
$l \in \widehat{I}(-\frac{1}{r})$ and
$l={\rm Max}\{1\leq k \leq L+1|n_k\neq 0\}
={\rm Min}\{1\leq k \leq L+1|m_k \neq 0\}$, we have
\begin{align}
&
\lim_{z_1 \to x^{i+j}z_2}
F_{i,j}^{(+)}\left(\frac{z_2}{z_1};a\right)
\Lambda_{m_1, m_2, \ldots,m_{L+1}}^{(i)}(z_1)
\Lambda_{n_1, n_2, \ldots, n_{L+1}}^{(j)}(z_2)\notag
\\
&=
-\frac{c(r,x) d_{m_l+n_l}(r,x)}{
d_{m_l}(r,x)d_{n_l}(r,x)}
\prod_{l=1}^{{\rm Min}(i,j)-1}\Delta_1(x^{2l+1})
\Lambda_{m_1+n_1, m_2+n_2, \ldots, m_{L+1}+n_{L+1}}^{(i+j)}(x^{i}z_2).
\label{appendixA:3}
\end{align}
$\bullet$~If $l$ satisfies 
$l \in \widehat{I}(-\frac{1}{r})$ and
$l={\rm Max}\{1\leq k \leq L+1|m_k\neq 0\}
={\rm Min}\{1\leq k \leq L+1|n_k \neq 0\}$, we have
\begin{align}
&
\lim_{z_1 \to x^{-(i+j)}z_2}
F_{i,j}^{(-)}\left(\frac{z_2}{z_1};a\right)
\Lambda_{m_1, m_2, \ldots,m_{L+1}}^{(i)}(z_1)
\Lambda_{n_1, n_2, \ldots, n_{L+1}}^{(j)}(z_2)\notag
\\
&=
\frac{c(r,x)d_{m_l+n_l}(r,x)}{d_{m_l}(r,x)d_{n_l}(r,x)}
\prod_{l=1}^{{\rm Min}(i,j)-1}\Delta_1(x^{2l+1})
\Lambda_{m_1+n_1, m_2+n_2, \ldots, m_{L+1}+n_{L+1}}^{(i+j)}(x^{-i}z_2).
\label{appendixA:4}
\end{align}
The remaining fusions vanish.

\section{Exchange relation}
\label{Appendix:exchange}
\setcounter{equation}{0}

In this appendix we give the exchange relations of $\Lambda_i(z)$
and $\Lambda_{m_1,m_2,\ldots,m_{L+1}}^{(i)}(z)$, 
which are obtained from Proposition \ref{prop:4-3}.
For $(m_1,m_2,\ldots,m_{L+1}) \in \hat{N}(\Phi)$ in
(\ref{def:N(Phi)}),
we set $i=m_1+m_2+\cdots+m_{L+1}$. We assume $i \geq 1$.
We calculate
\begin{eqnarray}
f_{1,i}\left(\frac{z_2}{z_1};a\right)\Lambda_l(z_1)\Lambda_{m_1,m_2,\ldots,m_{L+1}}^{(i)}(z_2)-
f_{i,1}\left(\frac{z_1}{z_2};a\right)
\Lambda_{m_1,m_2, \ldots, m_{L+1}}^{(i)}(z_2)\Lambda_l(z_1),
\label{appendixB:1}
\end{eqnarray}
where $a=D(0,L;\Phi)$.
\\
$\bullet$~If $l$ satisfies $m_l \neq 0$ and 
$l \in \widehat{I}(\frac{1-r}{r})$,
(\ref{appendixB:1}) is deformed as
\begin{eqnarray}
f_{1,i}\left(\frac{z_2}{z_1};a\right)\Lambda_l(z_1)\Lambda_{m_1,m_2,\ldots,m_{L+1}}^{(i)}(z_2)-
f_{i,1}\left(\frac{z_1}{z_2};a\right)
\Lambda_{m_1,m_2, \ldots,m_{L+1}}^{(i)}(z_2)\Lambda_l(z_1)=0.
\label{appendixB:2}
\end{eqnarray}
$\bullet$~If $l$ satisfies $m_l \neq 0$ and
$l \in \widehat{I}(-\frac{1}{r})$,
(\ref{appendixB:1}) is deformed as
\begin{align}
&\frac{c(r,x)d_{m_l+1}(r,x)}{d_1(r,x)d_{m_l}(r,x)}
:\Lambda_l(z_1)\Lambda_{m_1,m_2,\ldots,m_{L+1}}^{(i)}(z_2):
\notag
\\
&\times
\left(\delta\left(
x^{-i-1+2(m_1+m_2+\cdots+m_{l-1})}\frac{z_2}{z_1}\right)
-\delta
\left(x^{i+1-2(m_{l+1}+m_{l+2}+\cdots+m_{L+1})}
\frac{z_2}{z_1}\right)
\right).
\label{appendixB:3}
\end{align}
$\bullet$~If $l$ satisfies
$l < {\rm Min}\{1\leq k \leq L+1|m_k\neq 0\}$,
(\ref{appendixB:1}) is deformed as 
\begin{align}
c(r,x)
:\Lambda_l(z_1)\Lambda_{m_1,m_2,\ldots,m_{L+1}}^{(i)}(z_2):
\left(\delta\left(
\frac{x^{-i-1}z_2}{z_1}\right)
-\delta
\left(
\frac{x^{-i+1}z_2}{z_1}\right)
\right).
\label{appendixB:4}
\end{align}
$\bullet$~If $l$ satisfies
$l > {\rm Max}\{1\leq k \leq L+1|m_k\neq 0\}$,
(\ref{appendixB:1}) is deformed as
\begin{align}
c(r,x)
:\Lambda_l(z_1)\Lambda_{m_1,m_2,\ldots,m_{L+1}}^{(i)}(z_2):
\left(\delta\left(
\frac{x^{i-1}z_2}{z_1}\right)
-\delta
\left(
\frac{x^{i+1}z_2}{z_1}\right)
\right).
\label{appendixB:5}
\end{align}
$\bullet$~If $l$ satisfies $m_l=0$ and
${\rm Min}\{1\leq k \leq L+1|m_k\neq 0\}<l<
{\rm Max}\{1\leq k \leq L+1|m_k\neq 0\}$,
(\ref{appendixB:1}) is deformed as
\begin{align}
&c(r,x)
:\Lambda_l(z_1)\Lambda_{m_1,m_2,\ldots,m_{L+1}}^{(i)}(z_2):
\notag\\
&\times
\left(\delta\left(
x^{-i+1+2(m_1+m_2+\cdots+m_{l-2})}\frac{z_2}{z_1}\right)
-\delta
\left(
x^{i+1-2(m_l+m_{l+1}+\cdots+m_{L+1})}\frac{z_2}{z_1}\right)
\right).
\label{appendixB:6}
\end{align}

\end{appendix}
\end{document}